\documentclass[onecolumn]{autart}    

\usepackage{graphicx}
\usepackage{dsfont}
\usepackage{color}
\usepackage{amsfonts,amsmath}
\usepackage{caption}
\usepackage{subcaption}
\usepackage{hyperref}
\usepackage{float}
\usepackage{tikz,pgf}

\usepackage{pgfplots}

\newcommand{\lan}{\langle}
\newcommand{\ran}{\rangle}
\newcommand{\R}{\mathbb{R}}
\newcommand{\E}{\mathbb{E}}
\newcommand{\Law}{\mathcal{L}}

\newcommand{\1}{\mathds{1}}
\allowdisplaybreaks
\begin{document}

\begin{frontmatter}

\title{On modified Euler methods for McKean-Vlasov stochastic differential equations with super-linear coefficients\thanksref{footnoteinfo}} 

\thanks[footnoteinfo]{
Corresponding author: Qingshuo Song. All authors contributed equally.}

\author[UM]{Jiamin Jian}\ead{jiaminj@umich.edu},    
\author[WPI]{Qingshuo Song}\ead{qsong@wpi.edu},    
\author[CS]{Xiaojie Wang}\ead{x.j.wang7@csu.edu.cn},
\author[WPI]{Zhongqiang Zhang}\ead{zzhang7@wpi.edu}, 
\author[CS]{Yuying Zhao}\ead{zhaoyuying78@gmail.com} 

\address[UM]{Department of Mathematics, University of Michigan, Ann Arbor, MI 48109, USA}  

\address[WPI]{Department of Mathematical Sciences, Worcester Polytechnic Institute, Worcester, MA 01609, USA}  

\address[CS]{School of Mathematics and Statistics, HNP-LAMA, Central South University, Changsha, Hunan, P. R. China}

                                         
          
\begin{abstract}                          
We introduce a new class of numerical methods for solving McKean-Vlasov stochastic differential equations, which are relevant in the context of distribution-dependent or mean-field models, under super-linear growth conditions for both the drift and diffusion coefficients. Under certain non-globally Lipschitz conditions, the proposed numerical approaches have half-order convergence in the strong sense to the corresponding system of interacting particles associated with  McKean-Vlasov SDEs. By leveraging a result on the propagation of chaos, we establish the full convergence rate of the modified Euler approximations to the solution of the McKean-Vlasov SDEs. Numerical experiments are included to validate the theoretical results.
\end{abstract}

\end{frontmatter}

\section{Introduction}

McKean-Vlasov stochastic differential equations (MV-SDEs), also known as distribution-dependent or mean-field SDEs, extend traditional SDEs by incorporating the collective behavior of multiple interacting particles. 
Initially proposed by McKean \cite{mckean1966class}, \cite{mckean1967propagation}, this class of equations gained increasing attention following Dawson’s foundational work \cite{dawson1983critical} and the development of the Lion's derivative with respect to the measure variables \cite{lions2007cours}. 
Solving MV-SDEs is crucial in control theory as they model large-scale systems where individual components interact with the collective behavior of the group. They describe both the optimal path for mean-field controls \cite{bensoussan2013mean}, \cite{sznitman1991topics} and the equilibrium trajectory for mean-field games (see \cite{CD18I}, \cite{CD18II}, \cite{Mishura20}, \cite{wu2022stabilization} and references therein). This type of SDEs are essential for simplifying control problems in distributed systems, such as robots, power grids, or financial markets, by focusing on the statistical distribution of agents rather than individual interactions. They also account for uncertainty in noisy environments, enabling robust control strategies. Applications include multi-agent systems (see \cite{Benachour98}, \cite{Bossy97}) and other highly relevant fields like filtering (as highlighted in \cite{Crisan10}). By reducing system complexity and providing feedback control for distributed systems, MV-SDEs make large-scale optimization problems more tractable.

McKean–Vlasov SDEs are also widely used to model random phenomena across various scientific domains, including physics, biology, engineering, and neural activities, such as \cite{Baladron12}, \cite{Bolley11}, \cite{Bossy15}, \cite{CD18I}, \cite{CD18II}, \cite{Dreyer11}, \cite{Guhlke18}, \cite{mckean1966class} and references therein. As a result, there has been a notable surge in interest in related research.

Let $(\Omega, \mathcal{F}, \{ \mathcal{F}_t\}_{ t \geq 0 }, \mathbb{P})$ be a filtered probability space satisfying the usual conditions, where $\mathcal F_t$ is the augmented filtration of a standard $m$-dimensional Brownian motion $W = \{W(t)\}_{t \geq 0}$.
For a fixed terminal time $T > 0$, we consider the following McKean–Vlasov SDEs
\begin{equation}
\label{eq:MV_SDE}
X_t = X_0 + \int_0^t b \left(s, X_s,  \Law (X_s) \right) ds + \int_0^t \sigma \left(s, X_s,  \Law (X_s) \right) d W(s), \quad
t \in [0, T],
\ a.s.,
\end{equation}
where $\{\Law (X_t)\}_{t \geq 0}$ is the flow of deterministic marginal distributions of $X = \{X_t\}_{t \geq 0}$, $b:[0, T] \times \R^d \times \mathcal P_2(\R^d) \to \R^d$ denotes the drift function and $\sigma:[0, T] \times \R^d \times \mathcal P_2(\R^d) \to \R^{d \times m}$ is the diffusion function, expressed as $\sigma = (\sigma_1, \sigma_2, \dots, \sigma_m)$.
In this notation, $\sigma_j: [0, T] \times \R^d \times \mathcal P_2(\R^d) \to \R^{d}$ is the $j$-th column of $\sigma$.
Throughout this paper, the initial data
$X_0$ is a $\mathcal F_0$-measurable random variable in $\R^d$ independent of $W$.

In general, such equations rarely have explicit solutions available and one usually falls back on their numerical solutions. If the measure flow $\{\Law (X_t)\}_{t \geq 0}$ is known, then the coefficients $b$ and $\sigma$ are functions of time and space variables, and hence the MV-SDEs reduce to classical SDEs.
It is widely acknowledged that in scenarios where the coefficients of SDEs lack globally Lipschitz continuity and exhibit super-linear growth, the commonly used Euler-Maruyama numerical solution fails to attain finite moments, leading to divergence in both strong and weak senses. This issue has been well-documented in the literature, as evidenced by, e.g., \cite{Higham07}, \cite{Hutzen11}, \cite{Matting02},
\cite{KS17} and \cite{Mils05}. A similar divergence phenomenon, referred to as particle corruption, was observed in the context of MV-SDEs (see Section 4.1 in \cite{Reis22} for more details). 
Therefore, special care must be taken to construct and analyze convergent numerical approximations in a non-globally Lipschitz setting and recent years have witnessed a proper growth of the literature on this interesting topic \cite{Bao21},
\cite{liu2023tamed}, \cite{Reis22}, 
\cite{chen2022flexible}, \cite{gao2024convergence}, \cite{KN21}, \cite{Kumar22}, \cite{LMSWY23}, \cite{Liu23}, \cite{Neelima20}, \cite{chen2024euler}, \cite{chen2023wellposedness} and \cite{reisinger2022adaptive}.

Under local Lipschitz and linear growth conditions, the Euler numerical method for approximating MV-SDEs was analyzed in \cite{LMSWY23}. When the drift coefficients exhibit possible super-linear growth while the diffusion coefficients satisfy the linear growth condition, the moment boundedness and convergence rates of various numerical methods have been investigated in \cite{Bao21}, \cite{Reis22}, \cite{reisinger2022adaptive},\cite{fang2020adaptive}, \cite{KN21}, and 
\cite{liu2023tamed}. 
This analysis was further extended in \cite{Kumar22} to MV-SDEs with common noise and in \cite{Neelima20} to those with L\'evy processes, allowing the diffusion coefficients to also exhibit super-linear growth.
In addition, numerical methods have been proposed for solving a class of MV-SDEs with drift or diffusion components of convolution type. Specifically, \cite{chen2022flexible} and \cite{chen2024euler} addressed cases where both drifts and diffusions exhibit super-linear growth, while \cite{chen2023wellposedness} focused on drifts with super-linear growth and diffusion coefficients satisfying linear growth conditions. Finally, some works addressed numerical methods for McKean–Vlasov SDEs with H\"older continuous diffusion coefficient (see, e.g., \cite{Liu23}).
In terms of numerical methods, various approaches have been proposed in the literature:

\begin{itemize}
\item 
\textit{Explicit tamed Euler methods}, were introduced and studied in  \cite{Reis22}, \cite{liu2023tamed}, \cite{Liu23}, \cite{Neelima20}, \cite{Bao21}, \cite{KN21}, and \cite{Kumar22}, which rely on certain taming modifications of coefficients of MV-SDEs in a form such as $\frac{b(t, x, \mu)}{1+h^{\beta}|b(t, x, \mu)|} \, (0 < \beta \leq 1)$; 


\item 
\textit{Adaptive numerical methods} may serve as a viable alternative to tamed numerical solutions, particularly in numerically solving super-linear drift and diffusion coefficients, as demonstrated in \cite{reisinger2022adaptive};

\item 
\textit{Truncated method} was proposed in \cite{guo2024convergence} for the interacting particle system under a Khasminskii-type condition on the coefficients;

\item 
\textit{Projection-based particle method} was proposed in \cite{belomestny2018projected} to reduce the computational cost of solving MV-SDEs;

\item
\textit{Implicit numerical methods}, such as the backward Euler method \cite{Reis22} and 
\textit{split-step method} (see \cite{chen2022flexible}, \cite{chen2024euler}, \cite{chen2023wellposedness})
were utilized to approximate MV-SDEs with superlinear coefficients.
\end{itemize}

%
It is worthwhile to note that, when numerically approximating stable (dissipativity) systems, stability (or dissipativity) preserving methods are particularly vital. Usually, implicit methods have excellent stability properties and can preserve the dissipativity (long time stability) of the system, but at a price of expensive costs. 
A cheap option is to rely on explicit methods. 
However, as pointed out by \cite{chen2022flexible}, taming might destroy the strict dissipativity of the drift coefficients
and the usual tamed methods would be confronted with long time stability issues.
Our numerical experiments indicate different taming strategies would give different stability performances. 
Therefore, one should be careful with the choice of taming strategies for long-time simulations,
which, albeit interesting and important, turns out to be non-trivial.
Recently, several authors have made some attempts in this direction for usual SDEs.
By truncating monotonic functions, a recent work \cite{johnston2024strongly} proposed a polygonal (tamed) Euler method preserving the monotonicity of the drift coefficient. A new tamed Euler method was introduced by \cite{neufeld2025non-asymptotic} to long-time approximate invariant measures of the Langevin SDEs.

In this work, we focus on strong convergence analysis over finite time horizons for numerically solving MV-SDEs and leave the study of long time approximations for future work.
In the setting of the present article, the drift and diffusion coefficients are allowed to grow super-linearly in their spatial components. Rather than focusing on specific numerical solutions, we present a general framework to encompass a broader class of numerical methods. This allows us to establish moment boundedness and convergence rates within a general framework. A similar approach was adopted in \cite{lionnet2018convergence} for explicit numerical methods of forward-backward SDEs with drivers exhibiting polynomial growth.

The main contributions of this paper are summarized as follows:
\begin{itemize}
\item
\textit{New framework and new methods.}
We establish a new framework to admit novel numerical methods for MV-SDEs with super-linear drift and diffusion, such as the sin Euler method and tanh Euler method. As demonstrated by numerical experiments, the tanh Euler method has better stability properties than the other explicit tamed methods and always produces reliable approximations. 

\item \textit{Moment bound.}
Lemma \ref{l:moment_bound_euler_discrete_time} establishes the boundedness of moments of the newly proposed methods \eqref{eq:modified_euler_scheme} for MV-SDEs with super-linearly growing drift and diffusion coefficients. 
%

\item 
\textit{Convergence rate.}
We establish the strong convergence rates for a class of modified Euler methods in Theorem \ref{t:convergence_rate_euler}. It is shown that the proposed numerical methods have half-order convergence in the strong sense
of the corresponding system of interacting particles associated with the MV-SDEs. Moreover, the full convergence rate of the modified Euler approximation to the solution of the McKean-Vlasov SDEs \eqref{eq:MV_SDE} is provided in Corollary \ref{c:convergence_rate_euler} by leveraging a result on the propagation of chaos in Proposition \ref{l:propagation_of_chaos}.
\end{itemize}

The remainder of this article is structured as follows. In the forthcoming section, we present some necessary notations and 
our requirements on coefficients of MV-SDEs. A class of modified Euler methods and their uniform moment bounds are provided in Section \ref{sec:strong-con-MES}. Section \ref{s:convergence_rate} derives the strong convergence rate of the modified Euler approximations to the system of interacting particles. Finally, some numerical results are demonstrated in Section \ref{sec:numer-results}.


\section{Notations, assumptions, and preliminaries}
\label{s:Problem_setup}

In this section, we introduce notations and basic assumptions for the well-posedness of  MV-SDEs. The system of interacting particles and the corresponding result of propagation of chaos are also presented.

\subsection{Notations}

Let $|\cdot|$ and $\lan \cdot, \cdot\ran$ be the Euclidean norm and the inner product of vectors in $\R^d$, respectively. For a matrix $A$, we denote the Frobenius norm by $\|A\| = \sqrt{\text{tr}(A A^\top)}$, where $A^\top$ is the transpose of $A$ and $\text{tr}(\cdot)$ is the trace function of matrices. Let $\delta_x$ be the Dirac measure at a point $x \in \R^d$. 

To proceed, we denote $\mathcal P_2(\mathbb R^d)$ be the Wasserstein space of probability measures $\mu$ on $\mathbb R^d$ satisfying 
$\int_{\mathbb R^d} |x|^2 d \mu(x) < \infty$
endowed with $2$-Wasserstein metric  $\mathcal W_2(\cdot, \cdot)$ defined by
$$\mathcal W_2(\mu, \nu) = \inf_{\pi \in \Pi(\mu, \nu)} \left( \int_{\mathbb R^d} \int_{\mathbb R^d} |x - y|^2 d\pi(x, y) \right)^{\frac{1}{2}},$$
where $\Pi(\mu, \nu)$ is the collection of all probability measures on $\mathbb R^d \times \mathbb R^d$ with its marginals agreeing with $\mu$ and $\nu$.

\subsection{Well-posedness of MV-SDEs}

Next, we list the assumptions that are needed in this section. In the following, we use $L$ and $K$ to denote the generic constants which can be changed from line to line. Moreover, $p_0$ is denoted as a fixed positive constant that is sufficiently large and satisfies all the conditions specified in the inequalities presented in the theorem and lemmas of this paper.

\begin{assum}
\label{a:assumption_1_existence}
\begin{enumerate}
    \item[(A1)] $\mathbb E[|X_0|^{2 p_0}] < \infty$ for a fixed constant $p_0 > 1$.
    \vspace{4pt}
    
    \item [(A2)] There exists a constant $L > 0$ such that
    \begin{equation*}
    2 \lan x, b(t, x, \mu) \ran + (2 p_0 - 1) \|\sigma(t, x, \mu)\|^2 \leq L \left( 1 + |x|^2 + \mathcal W_2^2 (\mu, \delta_0) \right)
    \end{equation*}
    for all $t \in [0, T]$, $ x \in \R^d$ and $\mu \in \mathcal P_2(\R^d)$.
    \vspace{4pt}
    
    \item[(A3)] There exists a constant $L > 0$ such that
    \begin{equation*}
    2 \lan x - y, b(t, x, \mu) - b(t, y, \nu) \ran + \|\sigma(t, x, \mu) - \sigma(t, y, \nu) \|^2 \leq L \left(|x-y|^2 + \mathcal W_2^2 (\mu, \nu) \right)
    \end{equation*}
    for all $t \in [0, T] $, $x$, $y \in \R^d$ and $\mu$, $\nu \in \mathcal P_2(\R^d)$.
    \vspace{4pt}
    
    \item[(A4)] For every $t \in [0, T]$ and $\mu \in \mathcal P_2(\R^d)$, $b(t, \cdot, \mu)$ is a continuous function on $\R^d$ and for every $R > 0$ there exists $N_{R} \geq 0$ such that $\sup_{|x| \leq R} |b(t, x, \delta_0)| \leq N_{R}$ for all $t \in [0, T]$.
\end{enumerate}
\end{assum}

Assumption \ref{a:assumption_1_existence} plays a pivotal role in achieving existence, uniqueness, and moment boundedness concerning the McKean–Vlasov SDEs \eqref{eq:MV_SDE}. Detailed proof of the following result could be found in (\cite[Theorem 2.1]{Kumar22} with common noise).

\begin{prop}
(\cite[Theorem 2.1]{Kumar22})
\label{l:well_posedness_MV_sde}
Let assumptions (A1), (A2), (A3) and (A4) in Assumption \ref{a:assumption_1_existence} be satisfied. Then, there exists a unique solution to \eqref{eq:MV_SDE} and the following boundedness of moments hold
$$
\sup_{0 \leq t \leq T} \mathbb{E} \left[ \left| X_{t} \right|^{2p_{0}} \right] \leq K,
$$
where $p_0$ is from \textit{(A1)} and $K := K (L, \mathbb{E} [|X_{0}|^{2p_{0}}], d, m ) > 0$ is a constant. Moreover,
$$
\mathbb{E} \left[ \sup_{0 \leq t \leq T} \left| X_{t} \right|^{2q} \right] \leq K
$$
for all $q < p_0$.
\end{prop}

\subsection{The interacting particle system and propagation of chaos}

For a fixed $N \in \mathbb N$ and $i = 1, 2, \dots, N$, let $( W^{i}, X_0^{i})$ be $N$ independent copies of $(W, X_0)$. Note that, in the simulation of MV-SDEs, we need to approximate the measure $\Law (X_t)$ for all $t \geq 0$, which is not required in the case of classical SDEs. We consider the $N$-dimensional system of interacting particles
\begin{equation}
\label{eq:N_particle_system}
X_t^{i, N} = X_0^i + \int_0^t b \left(s, X_s^{i, N}, \mu_s^{X, N} \right) ds + \int_0^t \sigma \left(s, X_s^{i, N}, \mu_s^{X, N} \right) d W^{i} (s), \quad a.s.
\end{equation}
for all $t \in [0, T]$ and $i \in \{1, 2, \dots, N \}$, where $\mu_s^{X, N}$ is an empirical measure defined by
$$\mu_s^{X, N} (\cdot) = \frac{1}{N} \sum_{i = 1}^{N} \delta_{X_s^{i, N}} (\cdot).$$
Note that, $X_t^{i, N}$ in the $N$-dimensional system of interacting particles \eqref{eq:N_particle_system} is a proper approximation to $X_t$ in MV-SDEs \eqref{eq:MV_SDE} when $N$ is large enough. This result is called the propagation of chaos. Due to distribution dependence in \eqref{eq:MV_SDE}, we use the $N$-dimensional system of interacting particles \eqref{eq:N_particle_system} as a bridge to build the numerical approximations for MV-SDEs \eqref{eq:MV_SDE}. In order to present the propagation of chaos, we consider the following system of non-interacting particles:
\begin{equation}
\label{eq:non_interacting_particles}
X^{i}_t = X^{i}_0 + \int_0^t b \left(s, X^{i}_s,  \Law \left(X^{i}_s \right) \right) ds + \int_0^t \sigma \left(s, X^{i}_s,  \Law \left(X^{i}_s \right) \right) d W^{i}(s), \quad a.s.
\end{equation}
for all $t \in [0, T]$ and $i \in \{1, 2, \dots, N \}$. Note that, if the MV-SDEs \eqref{eq:MV_SDE} have a unique solution, then
$$\Law (X_t) = \Law \left(X^{i}_t \right), \quad \forall i = 1, 2, \dots, N, \quad \forall t \in [0, T].$$

Under Assumption \ref{a:assumption_1_existence}, Proposition 1 of \cite{Kumar22} asserts the result of the propagation of chaos.
\begin{prop}
\label{l:propagation_of_chaos}
(Propagation of chaos, \cite[Proposition 1]{Kumar22})
Let assumptions (A1), (A2), (A3) and (A4) in Assumption \ref{a:assumption_1_existence} hold with $p_0 > 2$. Then,
\begin{equation*}
\sup_{i \in \{1,2, \dots, N\}} \sup_{t \in [0, T]} \mathbb{E} \left[ \left\vert X_{t}^{i} - X_{t}^{i, N} \right\vert^{2} \right] \leq K 
\begin{cases} 
N^{-\frac{1}{2}}, & d < 4, \\
N^{-\frac{1}{2}} \ln (N), & d = 4, \\
N^{-\frac{2}{d},} & d > 4,
\end{cases}
\end{equation*}
where $K > 0$ is independent with $N$.
\end{prop}

\section{Modified Euler methods with moment bound}
\label{sec:strong-con-MES}

In this section, a class of modified Euler methods for the $N$-dimensional system of interacting particles \eqref{eq:N_particle_system} associated with the MV-SDEs \eqref{eq:MV_SDE} are proposed when the coefficients $b$ and $\sigma$ are allowed to grow super-linearly with respect to the state. Moreover, the boundedness of moments of the numerical approximation is also provided. 
\subsection{Modified Euler approximations}
\label{s:numerical_scheme}

Let $n \in \mathbb{N}$ be given, we construct a uniform mesh on $[0, T]$ with $h = T/n \in (0,1)$ being the stepsize and $t_k = k h$ for $k=0,1, \dots, n$. For $i = 1,2, \dots, N$ and $k = 0,1, \dots, n-1$, we consider modified Euler approximations
\begin{equation}
\label{eq:modified_euler_scheme}
\begin{aligned}
X_{t_{k+1}}^{i, N, n} &= X_{t_{k}}^{i, N, n} + \mathcal T_{1} \left( b \left(t_{k}, X_{t_{k}}^{i, N, n}, \mu_{t_{k}}^{X, N, n} \right), h \right) h \\
& \hspace{0.5in} + \sum_{r=1}^{m}  \mathcal T_{2} \left(\sigma_{r} \left(t_{k}, X_{t_{k}}^{i, N, n}, \mu_{t_{k}}^{X, N, n} \right), h \right) \Delta W_{r}^{i} \left(t_{k} \right)
\end{aligned}
\end{equation}
with $X_{t_0}^{i, N, n} = X_0^{i}$, where $\Delta W_{r}^{i}(t_{k}) = W_{r}^{i}(t_{k+1}) - W_{r}^{i}(t_{k})$ and $\mathcal T_{1}, \mathcal T_{2}$ are operators satisfying
$$
\mathcal T_{1}, \mathcal T_{2}: \mathbb{R}^{d} \times(0,1) \to \mathbb{R}^{d}.
$$
For the simplicity of the notation, we denote $X_{k}^{i, N, n} := X_{t_{k}}^{i, N, n}$ for all $i = 1,2, \dots, N$ and $k = 1, 2,  \dots, n$.


Next, we put some assumptions on the operators $\mathcal T_{1}$ and $\mathcal T_{2}$ in \eqref{eq:modified_euler_scheme}, in order to achieve the moment boundedness of the modified Euler approximations.
\begin{assum}
\label{a:euler_operator_t12}
\begin{enumerate}
    \item[(H1)] There exists a constant $L > 0$ such that
    \begin{equation*}
    \left\vert \mathcal T_{1}(x, h) \right\vert \leq \min \left\{  L h^{-2}, |x| \right\}, \quad
    \left\vert \mathcal T_{2}(x, h) \right\vert \leq \min \left\{ L h^{-\frac{3}{2}}, |x| \right\}
    \end{equation*}
    for all $x \in \mathbb{R}^{d}$ and $h \in (0,1)$.
    \vspace{4pt}
    
    \item[(H2)] There exist some constants $L, r_1, r_2 > 0$ such that
    \begin{equation*}
    \left\vert \mathcal T_{1}(x, h) - x \right\vert \leq L h^{r_1} |x|^{r_2}
    \end{equation*}
    for all $x \in \mathbb{R}^{d}$ and $h \in (0,1)$.
\end{enumerate}
\end{assum}
\begin{rem}
The condition \textit{(H1)} in Assumption \ref{a:euler_operator_t12} ensures that the maps 
$\mathcal{T}_1,\mathcal{T}_2$ are controlled by the linear growth and their values are also bounded by the inverse of the step size $h$. This is essentially used to aviod moment explosion and maintain stability even when the drift or diffusion terms exhibit polynomial growth. 
Moreover, the condition \textit{(H2)} serves as a consistency condition ensuring that the difference between $\mathcal T_1(x, h)$ and $x$ is sufficiently close in the sense that $\mathcal T_1(x, h) \to x$ as $h \to 0$ for fixed $x \in \mathbb R^d$.
To intuitively understand it, we depict the one-dimensional mappings $\mathcal{T}_i$, $i=1,2$ of different choices in Fig. \ref{fig:sketch-popu-exam-moment-bound}, satisfying Assumption \ref{a:euler_operator_t12}.
\end{rem}
\begin{center}
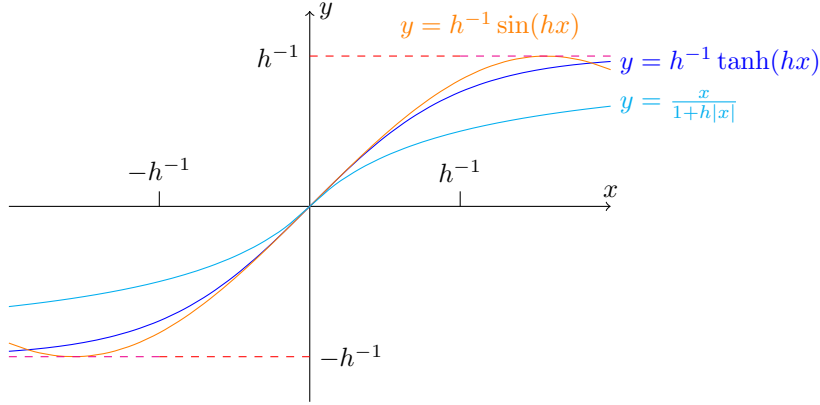

\begin{tikzpicture}[scale=0.20]
 \draw [->] (-20,0) -- (20,0) node [above] {$x$};
\draw [->] (0,-13) -- (0,13) node [right] {$y$};
\draw [color = blue, smooth , domain =-20:20]
  plot (\x, {10*(exp(\x/10)- exp(-\x/10))/(exp(\x/10)+ exp(-\x/10))})   node[right] {$y=h^{-1}\tanh(h x)$};
  
 \draw[color=orange, smooth, domain =-20:20]
 plot (\x, {10*sin(\x/10 r)})   node at (12,12) 
 {$y=h^{-1}\sin(h x)$};

  \draw [ color = cyan, smooth , domain =-20:20]
  plot (\x, {\x/(1 + (abs(\x)/10})   node[right] {$y=\frac{x}{1+h|x|}$};
  
  
  \draw [ color =red , dashed , domain =0:10]  plot (\x, 10) ;
  
  \draw [ color =magenta , dashed, domain =10:20]
  plot (\x, 10) ;

  \draw [ color =magenta, dashed, domain =-20:-10]
plot (\x, -10) ;

    \draw [ color =red, dashed , domain =-10:0]
  plot (\x, -10) ;

  \draw [solid] (-10,0.0) -- (-10,1)  node[above] {$-h^{-1}$};
  \draw [solid](10,0) -- (10,1) node[above] {$h^{-1}$};
  
  \draw (0,-10) node[right] {$-h^{-1}$};
  \draw (0,10) node[left] {$h^{-1}$};
\end{tikzpicture}
\captionof{figure}{Several choices for the operators $\mathcal{T}_i$, $i=1,2$.}
\label{fig:sketch-popu-exam-moment-bound}
\end{center}
Also, we mention that Assumption \ref{a:euler_operator_t12} just provides sufficient conditions used to derive the moment bounds of the numerical approximations. 
We present in Subsection \ref{s:example_mod_euler} some examples of modified Euler methods fulfilling Assumption \ref{a:euler_operator_t12}. In the literature, there are numerical methods that do not satisfy the condition \textit{(H1)} in Assumption \ref{a:euler_operator_t12}, whose moment boundedness can be derived in a different way (see Example \ref{ex:fully-tamed-Euler-scheme} below quoted from \cite{Neelima20}).

\subsection{Examples of modified Euler approximations}
\label{s:example_mod_euler}

The following are some examples of modified Euler type methods \eqref{eq:modified_euler_scheme}, where $\mathcal T_{1}$ and $\mathcal T_{2}$ are explicitly given.

\begin{exmp}[Drift-tamed Euler (DTE) \cite{Reis22}, \cite{Liu23}]
\label{ex:drift-tamed-Euler-scheme}
In \cite{Reis22}(\cite{Liu23}), the diffusion coefficient $\sigma$ is assumed to be globally Lipschitz continuous (or H\"older continuous). In their setting, the diffusion coefficient $\sigma$ does not need to be tamed, and a drift-tamed Euler method was introduced, where
\begin{equation}
\label{eq:numerical_approximation_DTE}
\mathcal T_{1}(x, h) = \frac{x}{1+h^{\lambda}|x| }, \quad \mathcal T_{2}(x, h) = x, \quad 0 < \lambda \leq \frac{1}{2}.
\end{equation} 
\end{exmp}

In this work, we propose three new modified Euler methods as follows.

\begin{exmp}[Modified Euler method (ME)] 
\label{ex:modified-Euler-scheme}
We propose a modified Euler method, where the operators $\mathcal T_{1}$ and $\mathcal T_{2}$ are given by
\begin{equation}
\label{eq:numerical_approximation_MES}
\mathcal T_{1}(x, h) = \frac{x}{1+h|x|^2 }, \quad  \mathcal T_{2}(x, h) = \frac{x}{1+h|x|^2 }.
\end{equation}
It is straightforward to verify that the conditions \textit{(H1)} and \textit{(H2)} in Assumption \ref{a:euler_operator_t12} are satisfied with $r_1 \in (0,1]$ and $r_2 = 3$. 
\end{exmp}

\begin{exmp}[tanh Euler method (TE)] 
\label{ex:tanh-Euler-scheme}
We introduce a $\tanh$ Euler method, where the operators $\mathcal T_{1}$ and $\mathcal T_{2}$ are as follows:
\begin{equation}
\label{eq:numerical_approximation_TES}
\mathcal T_{1}(x, h) = \frac{1}{h^{\alpha}} \tanh(h^{\alpha} x), \quad  \mathcal T_{2}(x, h) = \frac{1}{h^{\alpha}} \tanh(h^{\alpha} x),
\end{equation}
for some $\alpha \in (0, 3/2)$. The conditions \textit{(H1)} and \textit{(H2)} in Assumption \ref{a:euler_operator_t12} are satisfied when we choose $r_1 \in (0, \alpha]$ and $r_2 = 2$. 
\end{exmp}
In a similar way, we propose the sin Euler method as follows.
\begin{exmp}[sin Euler method (SE)] 
\label{ex:sin-Euler-scheme}
In the sin Euler method, the operators $\mathcal T_{1}$ and $\mathcal T_{2}$ are defined by
\begin{equation}
\label{eq:numerical_approximation_SES}
\mathcal T_{1}(x, h) = \frac{1}{h^{\alpha}} \sin(h^{\alpha} x), \quad \mathcal T_{2}(x, h) = \frac{1}{h^{\alpha}} \sin(h^{\alpha} x),
\end{equation}
for some $\alpha \in (0, 3/2)$. The conditions \textit{(H1)} and \textit{(H2)} in Assumption \ref{a:euler_operator_t12} are  fulfilled with $r_1 \in (0, \alpha)$ and $r_2 = 2$. 
\end{exmp}

Next, we provide tamed Euler type methods that do not satisfy Assumption \ref{a:euler_operator_t12}, but the moment boundedness and strong convergence rates of the numerical solutions have be obtained in the literature \cite{Neelima20}. 

\begin{exmp}[Fully-tamed Euler method (FTE) \cite{Neelima20}]
\label{ex:fully-tamed-Euler-scheme}
The author of \cite{Neelima20} proposed a drift-diffusion fully tamed Euler method, where 
$ \mathcal T_1, 
\mathcal T_2$ are given by
\begin{equation}
\label{eq:numerical_approximation_FTE}
\mathcal T_{1}( b(t, x, \mu), h) = 
\frac{ b(t, x, \mu)}{1+h^{\frac{1}{2}}|x|^{4\rho} }, \quad \mathcal T_{2}(\sigma_r(t, x, \mu), h) 
= \frac{ \sigma_r(t, x, \mu)}
{1+h^{\frac{1}{2}}|x|^{4 \rho} }.
\end{equation} 
Here $\rho$ comes from the growth condition (A6) of the drift $b$ below. It is not difficult to check that the mappings  $\mathcal T_1, 
\mathcal T_2$ do not obey \textit{(H1)} in Assumption \ref{a:euler_operator_t12}, but satisfy 
\begin{equation*}
\left\vert  \mathcal T_{1}( b(t, x, \mu), h) \right\vert \leq \min \big\{L h^{-\frac14}(1+|x|)+\mathcal W_2 (\mu, \delta_0), |b(t, x, \mu)| \big\}
\end{equation*} 
and
\begin{equation*}
\left\vert  \mathcal T_{2}( \sigma_r(t, x, \mu), h) \right\vert \leq \min \big\{  L h^{-\frac18}(1+|x|)+\mathcal W_2 (\mu, \delta_0), |\sigma_r(t, x, \mu)| \big\}.
\end{equation*} 
\end{exmp}
%

\subsection{Boundedness of moments of modified Euler approximations}
\label{sub:bound-moment-modi-method}

As demonstrated in \cite{Higham07}, \cite{Hutzen11}, \cite{Matting02}, \cite{Mils05}, it has been established that the numerical approximation generated by the Euler-Maruyama method lacks finite moments, which is of paramount significance for achieving convergence toward the desired system of interacting particles. Subsequently, we establish that the modified Euler numerical approximations, as defined by \eqref{eq:modified_euler_scheme}, possess bounded high-order moments.

In the following, we give some assumptions on the regularity of the coefficients in MV-SDEs to establish the moment boundedness  and the convergence rate of modified Euler approximations.

\begin{assum}
\label{a:euler_coefficient_b_sigma}
\begin{enumerate}
    \item[(A5)] For some $p_1 > 1$, there exists a constant $L > 0$ such that
    \begin{equation*}
    2 \lan x - y, b(t, x, \mu) - b(t, y, \nu) \ran + (2 p_1 - 1) \|\sigma(t, x, \mu) - \sigma(t, y, \nu) \|^2 \leq L \left(|x-y|^2 + \mathcal W_2^2 (\mu, \nu) \right)
    \end{equation*}
    for all $t \in [0, T], x, y \in \R^d$ and $\mu, \nu \in \mathcal P_2(\R^d)$.
    \vspace{4pt}
    
    \item[(A6)] There exist constants $L > 0$ and $\rho > 0$ such that
    \begin{equation*}
    |b(t, x, \mu) - b(t, y, \nu) | \leq L \left( \left( 1 + |x|^{2 \rho} + |y|^{2\rho} \right) |x - y|\right) + L \mathcal W_2 (\mu, \nu)
    \end{equation*}
    for all $t \in [0, T], x, y \in \R^d$ and $\mu, \nu \in \mathcal P_2(\R^d)$.
    \vspace{4pt}
    
    \item[(A7)] There exist a constant $L > 0$ such that
    \begin{equation*}
    |b(t, x, \mu) - b(s, x, \mu) | + \|\sigma(t, x, \mu) - \sigma(s, x, \mu) \| \leq L |t - s|^{\frac{1}{2}}
    \end{equation*}
    for all $t, s \in [0, T], x \in \R^d$ and $\mu \in \mathcal P_2(\R^d)$.
\end{enumerate}
\end{assum}

\begin{rem}
We mention that, assumption (A5) is stronger than the assumption (A3) since $2 p_1 - 1 > 1$. From assumptions (A5) and (A6), it follows that there exists a constant $K := K(L) > 0$
such that
\begin{equation*}
\|\sigma(t, x, \mu) - \sigma(t, y, v)\| \leq K  \left( \left(1+ |x|^{\rho} + |y|^{\rho} \right)|x-y| + \mathcal W_{2} (\mu, \nu) \right)
\end{equation*}
for all $t \in[0, T], x, y \in \mathbb{R}^{d}$, and $\mu, \nu \in \mathcal P_{2}(\mathbb{R}^{d})$.

Moreover, due to assumptions (A2), (A6) and (A7), there exists a constant $K:= K(L, T) > 0$ such that
\begin{equation}
\label{eq:growth_condition_b}
|b(t, x, \mu)| \leq K \left(1+|x|^{2 \rho+1} + \mathcal W_{2} \left(\mu, \delta_{0} \right) \right)
\end{equation}
and
\begin{equation}
\label{eq:growth_condition_sigma}
\|\sigma(t, x, \mu)\| \leq K \left( 1 + |x|^{\rho + 1} + \mathcal W_{2} \left(\mu, \delta_{0} \right) \right)
\end{equation}
for all $t \in [0, T], x \in \mathbb{R}^{d}$ and $\mu \in \mathcal  P_{2} (\mathbb{R}^{d})$. Note that, \eqref{eq:growth_condition_b} and \eqref{eq:growth_condition_sigma} provide the growth condition for the coefficients $b$ and $\sigma$, respectively.
\end{rem}

Also, we mention that, in the above settings, the drift and diffusion coefficients of the MV-SDEs are assumed to be $\mathcal{W}_2$-Lipschitz with respect to the measure component.
The following moment bound result, motivated by \cite{zhao2024one}, is essential for establishing the subsequent strong convergence result.
\begin{lem}
\label{l:moment_bound_euler_discrete_time}
Suppose assumptions (A1), (A2), (A6), (A7), (H1), (H2) hold. Then, for all $k = 0, 1, \dots, n$, there exists $\beta > 1$ and $K > 0$ independent of $n$  and $h$ such that
\begin{equation}
\label{eq:moment_boundness}
\sup_{i \in \{1, 2, \dots, N\}} \mathbb E \Big[ \Big\vert X_{k}^{i, N, n} \Big\vert^{2p} \Big] \leq K \big(1  + \mathbb E [ |X_0|^{2 p \beta}] \big), \quad \forall p \in \Big[1, \frac{2 \bar p - \mathcal G}{2 + 4 \mathcal G} \Big],
\end{equation}
where
\begin{equation}
\label{eq:constant_G}
\mathcal G := \mathcal G(\rho, r_1, r_2) = \max \Big\{6 \rho, \frac{(2 \rho + 1) r_2 - 1}{r_1} \Big\}
\end{equation}
with $\rho > 0$ from (A6) and $r_1, r_2 > 0$ from (H2), and $\bar{p}$ satisfies $p_0 \geq \bar{p} \geq 1 + \frac{5}{2} \mathcal G$.
\end{lem} 
The proof of this lemma can be found in Appendix \ref{appen:proof-moment-bound-modi-method}.
\begin{rem}
The boundedness of the moments of existing tamed numerical methods for MV-SDEs, such as those studied in \cite{Kumar22}, \cite{Neelima20}, was obtained based on an essential use of the following coercivity condition:
    \begin{equation}
    \label{eq:fra-mono-con-method}
    2 \lan x,  b_h(t, x, \mu)  \ran + (2 p_0 - 1) \big\| \sigma_h  (t, x, \mu)  \big\|^2 \leq L \left( 1 + |x|^2 + \mathcal W_2^2 (\mu, \delta_0) \right)
    \end{equation}
    for all $t \in [0, T], x \in \R^d$ and $\mu \in \mathcal P_2(\R^d)$. 
Here $b_h, \sigma_h$ are certain taming modifications of the drift and diffusion coefficients $b, \sigma$.
However, the modified Euler, tanh Euler and sin Euler methods proposed in Examples \ref{ex:modified-Euler-scheme}, \ref{ex:tanh-Euler-scheme}, \ref{ex:sin-Euler-scheme} fail to satisfy the condition \eqref{eq:fra-mono-con-method}. By formulating a different framework, here we employ new and different techniques to establish the desired moment bounds for these novel methods.
\end{rem}
%

\section{Strong convergence rate of modified Euler approximations}
\label{s:convergence_rate}

We prove the strong convergence rate of the modified Euler approximation \eqref{eq:modified_euler_scheme} in this section. Firstly, we provide a continuous-time version of modified Euler approximations for \eqref{eq:N_particle_system}. Let $\kappa_n(t) = \frac{\lfloor nt \rfloor}{n} = \sup \{s \in \{t_0, t_1, \dots, t_n\}, s \leq t\}$ for all $t \in [0, T]$. The modified Euler approximation in continuous time is given by
\begin{equation}
\label{eq:modified_euler_scheme_continuous_time}
\begin{aligned}
& X_{t}^{i, N, n} = X_0^i + \int_0^t \mathcal T_1 \left(b \left(\kappa_n(s), X_{\kappa_n(s)}^{i, N, n}, \mu_{\kappa_n(s)}^{X, N, n} \right), h \right) ds \\
& \hspace{0.8in} + \sum_{r = 1}^{m} \int_0^t \mathcal T_2 \left(\sigma_r \left(\kappa_n(s), X_{\kappa_n(s)}^{i, N, n}, \mu_{\kappa_n(s)}^{X, N, n} \right), h \right) dW_{r}^{i} (s)
\end{aligned}
\end{equation}
for all $t \in [0, T]$ and $i = 1, 2, \dots, N$.

Before proceeding with the proof of the rate of convergence of the modified Euler approximation \eqref{eq:modified_euler_scheme_continuous_time}, we establish some lemmas in what follows.

The following lemma provides an estimation between $X_t^{i, N, n}$ and $X_{\kappa_n(t)}^{i, N, n}$.

\begin{lem}
\label{l:difference_of_euler_numerical_scheme}
Under the same conditions of Lemma \ref{l:moment_bound_euler_discrete_time}, for all $i \in \{1,2, \dots, N\}$, $t \in[0, T]$ and $n$, $N \in \mathbb{N}$, we have the following inequality
$$
\mathbb{E} \Big[\left|X_{t}^{i, N, n} - X_{\kappa_n(t)}^{i, N, n} \right|^{2 p} \Big] \leq K h^{p}, \quad \forall p \in \Big[1, \frac{2 \bar p - \mathcal G}{(2 \rho + 1)(2 + 4 \mathcal G)} \Big],
$$
where $\mathcal G$ is given in \eqref{eq:constant_G} and $\bar{p}$ is a constant satisfying $p_0 \geq \bar{p} \geq (4 \rho + \frac{5}{2}) \mathcal G + 2 \rho + 1$. 
\end{lem}
\textit{Proof.}
Applying H\"older's inequality and the Burkholder-Davis-Gundy inequality, we have that for all $t \in [0, T]$,
\begin{equation*}
\begin{aligned}
& \mathbb{E}\left[\left|X_{t}^{i, N, n} - X_{\kappa_{n}(t)}^{i, N, n} \right|^{2 p} \right] \\
\leq \ & K \left(t - \kappa_{n}(t) \right)^{2 p - 1} \mathbb{E} \left[\int_{\kappa_{n}(t)}^{t} \left|\mathcal T_{1} \left(b \left(\kappa_{n}(s), X_{\kappa_{n}(s)}^{i, N, n}, \mu_{\kappa_{n}(s)}^{X, N, n}\right), h\right)\right|^{2 p} d s\right] \\
& \hspace{0.2in} + K \mathbb{E} \left[\left|\sum_{r=1}^{m} \int_{\kappa_{n}(t)}^{t} \mathcal T_{2}^{2} \left(\sigma_{r}\left(\kappa_{n}(s), X_{\kappa_{n}(s)}^{i, N, n}, \mu_{\kappa_{n}(s)}^{X, N, n} \right), h \right) d s \right|^{p} \right] \\
\leq \ & K h^{2 p - 1} \mathbb{E} \left[\int_{\kappa_{n}(t)}^{t} \left|\mathcal T_{1} \left( b \left(\kappa_{n}(s), X_{\kappa_{n}(s)}^{i, N, n}, \mu_{\kappa_{n}(s)}^{X, N, n}\right), h \right) \right|^{2 p} d s \right] \\
& \hspace{0.2in} + K h^{p - 1} \mathbb{E} \left[\sum_{r=1}^{m} \int_{\kappa_{n}(t)}^{t} \left|\mathcal T_{2} \left(\sigma_{r} \left(\kappa_{n}(s), X_{\kappa_{n}(s)}^{i, N, n}, \mu_{\kappa_{n}(s)}^{X, N, n} \right), h \right) \right|^{2 p} d s \right].
\end{aligned}
\end{equation*}
Under the assumption (\textit{H1}), and the growth condition for the coefficients $b$ and $\sigma$ in \eqref{eq:growth_condition_b} and \eqref{eq:growth_condition_sigma} respectively, it can be shown that
\begin{equation*}
\begin{aligned}
& \mathbb{E}\left[\left|X_{t}^{i, N, n} - X_{\kappa_{n}(t)}^{i, N, n} \right|^{2 p} \right] \\
\leq \ & K h^{2 p - 1} \mathbb{E} \left[\int_{\kappa_{n}(t)}^{t} \left| b \left(\kappa_{n}(s), X_{\kappa_{n}(s)}^{i, N, n}, \mu_{\kappa_{n}(s)}^{X, N, n} \right) \right|^{2 p} d s \right] \\
& \hspace{0.2in} + K h^{p - 1} \mathbb{E} \left[\sum_{r=1}^{m} \int_{\kappa_{n}(t)}^{t} \left|\sigma_{r} \left(\kappa_{n}(s), X_{\kappa_{n}(s)}^{i, N, n}, \mu_{\kappa_{n}(s)}^{X, N, n} \right) \right|^{2 p} d s \right] \\
\leq \ & K h^{2 p-1} \mathbb{E} \left[\int_{\kappa_{n}(t)}^{t} 1 + \left|X_{\kappa_{n}(s)}^{i, N, n}\right|^{2 p(2 \rho + 1)} + \left(\mathcal W_{2}^2 \left(\mu_{\kappa_{n}(s)}^{X, N, n}, \delta_{0} \right) \right)^{p} d s\right] \\
& \hspace{0.2in} + K h^{p - 1} \mathbb{E} \left[\int_{\kappa_{n}(t)}^{t} 1 + \left|X_{\kappa_{n}(s)}^{i, N, n} \right|^{2 p(\rho + 1)} + \left(\mathcal W_{2}^{2} \left(\mu_{\kappa_{n}(s)}^{X, N, n}, \delta_{0} \right) \right)^{p} d s \right].
\end{aligned}
\end{equation*}
From the identity \eqref{eq:w2_definition}, one can obtain
\begin{equation*}
\begin{aligned}
\mathbb{E}\left[\left|X_{t}^{i, N, n} - X_{\kappa_{n}(t)}^{i, N, n} \right|^{2 p} \right]
& \leq K h^{2 p}\left(1 + \sup_{s \in [\kappa_{n}(t), t]} \sup_{i \in \{1,2, \dots, N\}} \mathbb{E} \left[\left|X_{\kappa_{n}(s)}^{i, N, n} \right|^{2 p(2 \rho + 1)} \right] \right) \\
& \hspace{0.3in} + K h^{p} \left(1 + \sup_{s \in [\kappa_{n}(t), t]} \sup_{i \in \{1, 2, \dots, N\}} \mathbb{E} \left[ \left|X_{\kappa_{n}(s)}^{i, N, n} \right|^{2 p (\rho + 1)} \right] \right).
\end{aligned}
\end{equation*}
By the result of Lemma \ref{l:moment_bound_euler_discrete_time}, for all $i \in\{1,2, \dots, N\}$, $t \in[0, T]$, $n$, $N \in \mathbb N$ and $p \in [1, \frac{2 \bar p - \mathcal G}{(2 \rho + 1)(2 + 4 \mathcal G)}]$, we have
\begin{equation*}
\mathbb{E}\left[\left|X_{t}^{i, N, n} - X_{\kappa_{n}(t)}^{i, N, n} \right|^{2 p} \right] \leq K h^{p} \left(1 + \mathbb{E} \left[ \left|X_{0}\right|^{2 \beta p (2 \rho + 1)} \right] \right) \leq K h^{p},
\end{equation*}
as required.
\qed

The following lemma gives the boundedness of moments for the modified Euler approximation \eqref{eq:modified_euler_scheme_continuous_time}.

\begin{lem}
\label{l:moment_bound_euler_continuous_time}
Under the same conditions of Lemma \ref{l:moment_bound_euler_discrete_time}, there exist $\beta > 1$ and $K > 0$ independent of $n$ and $h$ such that
\begin{equation*}
\sup_{i \in \{1, 2, \dots, N\}} \sup_{t \in [0, T]} \mathbb E \left[ \left\vert X_{t}^{i, N, n} \right\vert^{2p} \right] \leq K \left(1  + \mathbb E \left[ |X_0|^{2 p \beta}  \right] \right), \quad \forall p \in \left[1, \frac{2 \bar p - \mathcal G}{(2 \rho + 1)(2 + 4 \mathcal G)} \right],
\end{equation*}
where $\mathcal G$ is given in \eqref{eq:constant_G} and $\bar{p}$ is a constant satisfying $p_0 \geq \bar{p} \geq (4 \rho + \frac{5}{2}) \mathcal G + 2 \rho + 1$.
\end{lem}
\textit{Proof.}
Applying Lemma \ref{l:moment_bound_euler_discrete_time} and Lemma \ref{l:difference_of_euler_numerical_scheme}, we obtain the desired result that
\begin{equation*}
\begin{aligned}
\sup_{i \in \{1, 2, \dots, N\}} \sup_{t \in [0, T]} \mathbb{E} \left[\left|X_{t}^{i, N, n} \right|^{2 p} \right]
& \leq \sup_{i \in \{1, 2, \dots, N\}} \sup_{t \in [0, T]} K \mathbb{E} \left[\left|X_{t}^{i, N, n} - X_{\kappa_{n}(t)}^{i, N, n} \right|^{2 p}\right] \\
& \hspace{0.5in} + \sup_{i \in \{1, 2, \dots, N\}} \sup_{t \in [0, T]} K \mathbb{E} \left[\left|X_{\kappa_{n} (t)}^{i, N, n} \right|^{2 p} \right] \\
& \leq K \left(1 + \mathbb{E} \left[\left|X_{0} \right|^{2 p \beta} \right] \right).
\end{aligned}
\end{equation*}
\qed

Next, we provide the assumptions for the proof of the convergence rate of the time-continuous approximation \eqref{eq:modified_euler_scheme_continuous_time}.

\begin{assum}
\label{a:euler_operator_t12_convergence}
\begin{enumerate}
\item[(H3)] There exist some constants $L > 0, r_2 > 0$ and $r_1, r_3 \geq \frac{1}{2}$ such that
\begin{equation*}
\left\vert \mathcal T_{1}(x, h) - x \right\vert \leq L h^{r_1} |x|^{r_2}, \quad
\left\vert \mathcal T_{2}(x, h) - x \right\vert \leq L h^{r_3} |x|^{r_2}
\end{equation*}
for all $x \in \mathbb{R}^{d}$ and $h \in (0,1)$.
\end{enumerate}
\end{assum}

Note that, in the assumption (\textit{H2}) of Assumption \ref{a:euler_operator_t12}, we only need $r_1> 0$. But in assumption (\textit{H3}) of Assumption \ref{a:euler_operator_t12_convergence}, we require that $r_1, r_3 \geq \frac{1}{2}$. Thus, the assumption \textit{(H3)} is slightly stronger than \textit{(H2)} as $h \in (0, 1)$.


\begin{lem}
\label{l:estimation_difference}
Suppose the assumptions (A1), (A2), (A6), (A7), (H1), (H3) hold. Then, for all $p \in [1,\tilde{p}]$, there exists $K > 0$ independent of $n$ and $h$ such that
\begin{equation*}
\mathbb{E} \left[ \left| b \left(t, X_{t}^{i, N, n}, \mu_{t}^{X, N, n} \right) - \mathcal T_{1} \left( b \left(\kappa_{n}(t), X_{\kappa_{n}(t)}^{i, N, n}, \mu_{\kappa_{n}(t)}^{X, N, n} \right), h \right) \right|^{2 p} \right] \leq K h^{p}
\end{equation*}
and
\begin{equation*}
\sup_{r\in \{ 1, 2, \dots, m\}} \mathbb{E} \left[ \left|\sigma_{r} \left(t, X_{t}^{i, N, n}, \mu_{t}^{X, N, n}\right) - \mathcal T_{2}\left( \sigma_{r} \left(\kappa_{n}(t), X_{\kappa_{n}(t)}^{i, N, n}, \mu_{\kappa_{n}(t)}^{X, N, n} \right), h \right) \right|^{2 p} \right] \leq K h^{p},
\end{equation*} 
where 
\begin{equation}
\label{eq:constant_tilde_p}
\tilde{p} = \min \left\{\frac{2 \bar p - \mathcal G}{r_2 (2 \rho + 1) (2 + 4 \mathcal G)}, \, \frac{2 \bar p - \mathcal G}{4 \rho (2 \rho + 1)(2 + 4 \mathcal G)}, \, \frac{2 \bar p - \mathcal G}{(2 \rho + 1)(2 + 4 \mathcal G)} \right\}
\end{equation}
with $\mathcal G$ is given in \eqref{eq:constant_G} and $\bar p \leq p_0$ is a constant such that
\begin{footnotesize}
\begin{equation*}
\bar p \geq \max \left\{\left(16 \rho^2 + 8 \rho + \frac{1}{2} \right) \mathcal G + 8 \rho^2 + 4 \rho, \, \left(2 r_2 (2 \rho + 1) + \frac{1}{2} \right) \mathcal G + r_2 (2 \rho + 1), \,  \left(4 \rho + \frac{5}{2} \right) \mathcal G + 2 \rho + 1 \right\}.
\end{equation*}
\end{footnotesize}
\end{lem}

\textit{Proof.}
By Lemma \ref{l:difference_of_euler_numerical_scheme}, H\"older's inequality, (\textit{A6}), (\textit{A7}), (\textit{H3}), and the growth condition \eqref{eq:growth_condition_b}, we have 
\begin{equation*}
\begin{aligned}
& \mathbb{E} \left[ \left| b \left(t, X_{t}^{i, N, n}, \mu_{t}^{X, N, n} \right) - \mathcal T_{1} \left( b \left(\kappa_{n}(t), X_{\kappa_{n}(t)}^{i, N, n}, \mu_{\kappa_{n}(t)}^{X, N, n}\right), h \right) \right|^{2 p} \right] \\
\leq \ & K \mathbb{E} \left[ \left| b \left(t, X_{t}^{i, N, n}, \mu_{t}^{X, N, n} \right) - b \left(t, X_{\kappa_{n}(t)}^{i, N, n}, \mu_{\kappa_{n}(t)}^{X, N, n} \right) \right|^{2 p} \right. \\
& \hspace{0.3in} + \left| b \left(t, X_{\kappa_{n}(t)}^{i, N, n}, \mu_{\kappa_{n}(t)}^{X, N, n} \right) - b \left(\kappa_{n}(t), X_{\kappa_{n}(t)}^{i, N, n}, \mu_{\kappa_{n}(t)}^{X, N, n} \right) \right|^{2 p} \\
& \hspace{0.3in} + \left. \left| b \left(\kappa_{n}(t), X_{\kappa_{n}(t)}^{i, N, n}, \mu_{\kappa_{n}(t)}^{X, N, n} \right) - \mathcal T_{1} \left( b \left(\kappa_{n}(t), X_{\kappa_{n}(t)}^{i, N, n}, \mu_{\kappa_{n}(t)}^{X, N, n} \right), h \right) \right|^{2 p} \right] \\
\leq \ & K \mathbb{E} \bigg[ \left(1 + \left|X_{t}^{i, N, n} \right|^{2 \rho} + \left|X_{\kappa_{n}(t)}^{i, N, n} \right|^{2 \rho} \right)^{2 p} \left|X_{t}^{i, N, n} - X_{\kappa_{n}(t)}^{i, N, n} \right|^{2p} \\
& \hspace{0.3in} + \mathcal W_{2}^{2p} \left(\mu_{t}^{X, N, n}, \mu_{\kappa_{n}(t)}^{X, N, n} \right) + \left|t - \kappa_{n}(t) \right|^{p} + h^{2 p r_{1}} \left|b \left(\kappa_{n}(t), X_{\kappa_{n}(t)}^{i, N, n}, \mu_{\kappa_{n}(t)}^{X, N, n} \right) \right|^{2 p r_{2}}\bigg] \\
\leq \ & K h^p \mathbb{E} \left[\left(1 + \left|X_{t}^{i, N, n} \right|^{8p \rho} + \left|X_{\kappa_{n}(t)}^{i, N, n} \right|^{8p \rho} \right)  \right] + K \sup_{i \in \{1, 2, \dots, N\}} \mathbb{E} \left[ \left| X_{t}^{i, N, n} - X_{\kappa_{n}(t)}^{i, N, n} \right|^{2 p} \right] + K h^{p} \\
& \hspace{0.3in} + K h^{2 p r_{1}} \mathbb{E} \left[1 + \left|X_{\kappa_{n}(t)}^{i, N, n} \right|^{2 p r_{2}(2 \rho + 1)} + \mathcal W_{2}^{2 p r_2} \left(\mu_{\kappa_{n}(s)}^{X, N, n}, \delta_{0} \right) \right]
\end{aligned}
\end{equation*}
since
\begin{equation*}
\mathcal W_{2}^{2} \left(\mu_{t}^{X, N, n}, \mu_{\kappa_{n}(t)}^{X, N, n} \right) \leq \frac{1}{N} \sum_{i=1}^{N}\left|X_{t}^{i, N, n} - X_{\kappa_{n}(t)}^{i, N, n} \right|^{2} \leq \sup_{i \in \{1, 2, \dots, N\}} \left|X_{t}^{i, N, n} - X_{\kappa_{n}(t)}^{i, N, n} \right|^{2}.
\end{equation*}
Then, Lemmas \ref{l:moment_bound_euler_discrete_time}, \ref{l:difference_of_euler_numerical_scheme} and \ref{l:moment_bound_euler_continuous_time} together show that
\begin{equation*}
\begin{aligned}
& \mathbb{E} \left[ \left| b \left(t, X_{t}^{i, N, n}, \mu_{t}^{X, N, n} \right) - \mathcal T_{1} \left( b \left(\kappa_{n}(t), X_{\kappa_{n}(t)}^{i, N, n}, \mu_{\kappa_{n}(t)}^{X, N, n}\right), h \right) \right|^{2 p} \right] \\
\leq \ &  K h^{p} \left(1 + \mathbb{E} \left[ \left|X_{0} \right|^{8 p \rho \beta} \right] \right) + K h^{p} + K h^{2 p r_{1}} \left(1 + \mathbb{E} \left[\left|X_{0} \right|^{2 p \beta r_{2}(2 \rho + 1)} \right] \right) \\
& \hspace{0.3in} + K h^{2 p r_{1}} \mathbb{E} \left[\left(\frac{1}{N} \sum_{i=1}^{N} \left|X_{\kappa_{n}(s)}^{i, N, n} \right|^{2} \right)^{p r_{2}} \right] \\
\leq \ & K h^{p}
\end{aligned}
\end{equation*}
since $r_1 \geq \frac{1}{2}$. The proof is completed by performing a similar calculation for $\sigma_r$ for all $r = 1, 2, \dots, m$ with $r_2 \geq \frac{1}{2}$.
\qed

Now, we are ready to prove the strong convergence rate of order $\frac{1}{2}$ for the modified Euler approximation \eqref{eq:modified_euler_scheme_continuous_time} in the $L^p$ sense. 

\begin{thm}
\label{t:convergence_rate_euler}
Suppose assumptions (A1), (A2), (A5), (A6), (A7), (H1), (H3) are satisfied. Then, the modified Euler approximation \eqref{eq:modified_euler_scheme_continuous_time} converges to the solution of the interacting particle system \eqref{eq:N_particle_system} in a strong sense with $L^p$ convergence rate given by
$$
\sup_{i \in \{1, 2, \dots, N\}} \sup_{t \in [0, T]} \mathbb{E} \left[ \left|X_{t}^{i, N} -X_{t}^{i, N, n} \right|^{2 p} \right] \leq K h^{p}
$$
for all $1 \leq p \leq \min \{\frac{1}{2} p_{1} + \frac{1}{4}, \tilde{p} \}$, where $p_1$ comes from assumption \textit{(A5)}, $\tilde{p}$ is given by \eqref{eq:constant_tilde_p} and $K > 0$ does not depend on $n, N \in \mathbb N$.
\end{thm}
\textit{Proof.}
From \eqref{eq:N_particle_system}, \eqref{eq:modified_euler_scheme_continuous_time}, and Ito's formula, it follows that
\begin{equation*}
\begin{aligned}
& \left|X_{t}^{i, N} - X_{t}^{i, N, n} \right|^{2 p} \\
= \ & 2 p \int_{0}^{t} \left|X_{s}^{i, N} - X_{s}^{i, N, n} \right|^{2 p - 2} \left\langle X_{s}^{i, N} - X_{s}^{i, N, n}, b \left(s, X_{s}^{i, N}, \mu_{s}^{X, N} \right) - \mathcal T_{1} \left(b \left(\kappa_{n}(s), X_{\kappa_{n}(s)}^{i, N, n}, \mu_{\kappa_{n}(s)}^{X, N, n} \right), h \right) \right\rangle d s \\
& \hspace{0.1in} + 2 p \int_{0}^{t} \left|X_{s}^{i, N} - X_{s}^{i, N, n} \right|^{2 p - 2} \sum_{r=1}^{m} \left\langle X_{s}^{i, N} - X_{s}^{i, N, n}, \sigma_{r} \left(s, X_{s}^{i, N}, \mu_{s}^{X, N} \right) - \mathcal T_{2} \left(\sigma_{r} \left(\kappa_{n}(s), X_{\kappa_{n}(s)}^{i, N, n}, \mu_{\kappa_{n}(s)}^{X, N, n} \right) , h \right)  \right\rangle d W_{r}^{i}(s) \\
& \hspace{0.1in} +  p (2 p - 1) \int_{0}^{t} \left|X_{s}^{i, N} - X_{s}^{i, N, n} \right|^{2 p-2} \sum_{r=1}^{m}\left|\sigma_{r} \left(s, X_{s}^{i, N}, \mu_{s}^{X, N} \right) - \mathcal T_{2} \left(\sigma_{r} \left(\kappa_{n}(s), X_{\kappa_{n}(s)}^{i, N, n}, \mu_{\kappa_{n}(s)}^{X, N, n} \right) , h \right) \right|^{2} d s.
\end{aligned}
\end{equation*}
Taking expectation and applying the Cauchy–Schwarz inequality, we have
\begin{equation*}
\begin{aligned}
& \mathbb{E} \left[\left|X_{t}^{i, N} - X_{t}^{i, N, n} \right|^{2 p} \right] \\ 
\leq \ & p \mathbb{E} \left[\int_{0}^{t} \left|X_{s}^{i, N} - X_{s}^{i, N, n} \right|^{2p - 2} \Big\{ \Big\langle 2 \left(X_{s}^{i, N} - X_{s}^{i, N, n} \right), b \left(s, X_{s}^{i, N}, \mu_{s}^{X, N} \right) - b \left(s, X_{s}^{i, N, n}, \mu_{s}^{X, N, n} \right) \right.  \\
& \hspace{0.2in} + b \left(s, X_{s}^{i, N, n}, \mu_{s}^{X, N, n} \right) - \mathcal T_{1} \left(b \left(\kappa_{n}(s), X_{\kappa_{n}(s)}^{i, N, n}, \mu_{\kappa_{n}(s)}^{X, N, n} \right), h \right)\Big\rangle \\
& \hspace{0.2in} + (2 p - 1) \sum_{r=1}^{m} \Big| \sigma_{r} \left(s, X_{s}^{i, N}, \mu_{s}^{X, N} \right) - \sigma_{r} \left(s, X_{s}^{i, N, n}, \mu_{s}^{X, N, n} \right) \\
& \hspace{0.2in} \left. \left. + \sigma_{r} \left(s, X_{s}^{i, N, n}, \mu_{s}^{X, N, n} \right) - \mathcal T_{2} \left(\sigma_{r} \left(\kappa_{n}(s), X_{\kappa_{n}(s)}^{i, N, n}, \mu_{\kappa_{n}(s)}^{X, N, n} \right) , h \right) \Big|^2 \right. \Big\} ds \right].
\end{aligned}
\end{equation*}
Noting $p \leq \frac{1}{2} p_1 + \frac{1}{4}$, we derive that
\begin{equation*}
\begin{aligned}
& \mathbb{E} \left[\left|X_{t}^{i, N} - X_{t}^{i, N, n} \right|^{2 p} \right] \\
\leq \ & p \mathbb{E} \bigg[\int_{0}^{t} \left|X_{s}^{i, N} - X_{s}^{i, N, n} \right|^{2p - 2} \Big\{ \left\langle 2 \left(X_{s}^{i, N} - X_{s}^{i, N, n} \right), b \left(s, X_{s}^{i, N}, \mu_{s}^{X, N} \right) - b \left(s, X_{s}^{i, N, n}, \mu_{s}^{X, N, n} \right) \right\rangle \\
& \hspace{0.1in} + (2 p_1 - 1) \sum_{r=1}^{m} \left| \sigma_{r} \left(s, X_{s}^{i, N}, \mu_{s}^{X, N} \right) - \sigma_{r} \left(s, X_{s}^{i, N, n}, \mu_{s}^{X, N, n} \right) \right|^2 \Big\} ds \bigg] \\
& \hspace{0.1in} + K \mathbb{E} \left[\int_{0}^{t} \left|X_{s}^{i, N} - X_{s}^{i, N, n} \right|^{2p - 2} \left\langle \left(X_{s}^{i, N} - X_{s}^{i, N, n} \right), b \left(s, X_{s}^{i, N, n}, \mu_{s}^{X, N, n} \right) - \mathcal T_{1} \left(b \left(\kappa_{n}(s), X_{\kappa_{n}(s)}^{i, N, n}, \mu_{\kappa_{n}(s)}^{X, N, n} \right), h \right) \right\rangle ds \right] \\
& \hspace{0.1in} + K \mathbb{E} \bigg[\int_{0}^{t} \left|X_{s}^{i, N} - X_{s}^{i, N, n} \right|^{2p - 2} \sum_{r=1}^{m} \left| \sigma_{r} \left(s, X_{s}^{i, N, n}, \mu_{s}^{X, N, n} \right) - \mathcal T_{2} \left(\sigma_{r} \left(\kappa_{n}(s), X_{\kappa_{n}(s)}^{i, N, n}, \mu_{\kappa_{n}(s)}^{X, N, n} \right) , h \right) \right|^2  ds \bigg]
\end{aligned}
\end{equation*}
for all $t \in[0, T], i \in \{1,2, \dots, N\}$ and $n, N \in \mathbb N$. By using the Cauchy–Schwarz inequality, Young's inequality, and Assumption (\textit{A5}), one can obtain
\begin{equation*}
\begin{aligned}
& \mathbb{E} \left[\left|X_{t}^{i, N} - X_{t}^{i, N, n} \right|^{2 p} \right] \\
\leq \ &  K \mathbb{E} \left[ \int_{0}^{t} \left|X_{s}^{i, N} - X_{s}^{i, N, n} \right|^{2 p} d s \right] + K \mathbb{E} \left[\int_{0}^{t} \mathcal W_{2}^{2 p} \left(\mu_{s}^{X, N}, \mu_{s}^{X, N, n} \right) d s \right] \\
& \hspace{0.2in} + K \mathbb{E} \left[ \int_{0}^{t} \left| b \left(s, X_{s}^{i, N, n}, \mu_{s}^{X, N, n} \right) - \mathcal T_{1} \left(b \left(\kappa_{n}(s), X_{\kappa_{n}(s)}^{i, N, n}, \mu_{\kappa_{n}(s)}^{X, N, n} \right), h \right) \right|^{2p} d s \right] \\
& \hspace{0.2in} + K \mathbb{E} \bigg[ \int_{0}^{t} \sum_{r=1}^{m} \left| \sigma_{r} \left(s, X_{s}^{i, N, n}, \mu_{s}^{X, N, n} \right) - \mathcal T_{2} \left(\sigma_{r} \left(\kappa_{n}(s), X_{\kappa_{n}(s)}^{i, N, n}, \mu_{\kappa_{n}(s)}^{X, N, n} \right) , h \right) \right|^{2p}  ds \bigg].
\end{aligned}
\end{equation*}
Lemma \ref{l:estimation_difference} along with the following estimate
$$\mathcal W_{2}^{2} \left(\mu_{s}^{X, N}, \mu_{s}^{X, N, n} \right) \leq \frac{1}{N} \sum_{i=1}^{N}\left|X_{s}^{i, N} - X_{s}^{i, N, n} \right|^{2}$$
for all $s \in [0, t]$ yields 
\begin{equation*}
\begin{aligned}
& \mathbb{E} \left[\left|X_{t}^{i, N} - X_{t}^{i, N, n} \right|^{2 p} \right] \\
\leq \ & K \mathbb{E} \left[ \int_{0}^{t} \left|X_{s}^{i, N} - X_{s}^{i, N, n} \right|^{2 p} d s \right] + K \mathbb{E} \left[\int_{0}^{t} \left(\frac{1}{N} \sum_{i=1}^{N}\left|X_{s}^{i, N} - X_{s}^{i, N, n} \right|^{2} \right)^{p} d s \right] + K \mathbb{E} \left[\int_{0}^{t} h^p  ds \right] \\
\leq \ & K \mathbb{E} \left[ \int_{0}^{t} \sup_{i \in \{1, 2, \dots, N\} }  \left|X_{s}^{i, N} - X_{s}^{i, N, n} \right|^{2 p} d s \right] + K h^p.
\end{aligned}
\end{equation*}
Therefore, the following inequality
\begin{equation*}
\sup_{i \in \{1, 2, \dots, N\} }  \sup_{r \in [0, t]} \mathbb{E} \left[\left|X_{r}^{i, N} - X_{r}^{i, N, n} \right|^{2 p} \right] \leq  K \int_{0}^{t} \sup_{i \in \{1, 2, \dots, N\} } \sup_{r \in [0, s]} \mathbb{E} \left[ \left|X_{r}^{i, N} - X_{r}^{i, N, n} \right|^{2 p} \right] d s + K h^p
\end{equation*}
holds for all $t \in [0, T]$, $i \in \{1, 2, \dots, N\}$ and $n, N \in \mathbb N$. The application of Gr\"onwall's inequality implies
\begin{equation*}
\sup_{i \in \{1, 2, \dots, N\} }  \sup_{t \in [0, T]} \mathbb{E} \left[\left|X_{t}^{i, N} - X_{t}^{i, N, n} \right|^{2 p} \right] \leq K h^p.
\end{equation*}
The proof is thus finished.
\qed

\begin{dis}
\label{dis:para-value-example}
Under the assumptions of Theorem \ref{t:convergence_rate_euler}, the strong convergence rates of the modified Euler approximations \eqref{eq:modified_euler_scheme_continuous_time} can be derived, depending on the specific modifications made to the coefficients. Specifically:
\begin{itemize}
    \item For Example \ref{ex:modified-Euler-scheme}, we confirm that assumption (H3) is satisfied with $r_1 =\frac12 $, $r_2=2$, $r_3=\frac12$.
     \item For Example \ref{ex:tanh-Euler-scheme}, we demonstrate that assumption (H3) is fulfilled with $r_1 =\frac12 $, $r_2=2$, $r_3=\frac12$.
     \item For Example \ref{ex:sin-Euler-scheme}, it is clear that assumption (H3) holds with $r_1 =\frac12 $, $r_2=2$, $r_3=\frac12$.
\end{itemize}
As a consequence of Theorem \ref{t:convergence_rate_euler}, the strong convergence rate of order $1/2$ is achieved for numerical methods defined by Examples \ref{ex:modified-Euler-scheme}, \ref{ex:tanh-Euler-scheme} and \ref{ex:sin-Euler-scheme}.   
\end{dis}

To ensure the completeness of our work, we now are ready to present the full numerical approximation errors of the modified Euler approximation \eqref{eq:modified_euler_scheme_continuous_time} to the solution of McKean-Vlasov SDE \eqref{eq:MV_SDE}. It is a straightforward result from the combination of Theorem \ref{t:convergence_rate_euler} and Proposition \ref{l:propagation_of_chaos}. 
\begin{cor}
\label{c:convergence_rate_euler}
Suppose assumptions (A1), (A2), (A3), (A4), (A5), (A6), (A7), (H1), (H3) are satisfied with $p_0 > 2$. Then, the modified Euler approach \eqref{eq:modified_euler_scheme_continuous_time} converges to the solution of McKean-Vlasov SDEs \eqref{eq:MV_SDE} in a strong sense with $L^2$ convergence rate given by
\begin{equation*}
\sup_{i \in \{1, 2, \dots, N\}} \sup_{t \in [0, T]} \mathbb{E} \left[ \left|X_{t}^{i} -X_{t}^{i, N, n} \right|^{2} \right] \leq K 
\begin{cases} 
h + N^{-\frac{1}{2}}, & d < 4, \\
h + N^{-\frac{1}{2}} \ln (N), & d = 4, \\
h + N^{-\frac{2}{d},} & d > 4,
\end{cases}
\end{equation*}
where $K > 0$ is independent of $n$ and $N \in \mathbb N$.
\end{cor}


\section{Numerical experiments}
\label{sec:numer-results} 
Some numerical experiments are conducted to illustrate the previous theoretical findings. 
Newton's method is used to solve implicit algebraic equations if necessary.

We begin by considering two examples to illustrate the mean square convergence rates of the modified Euler method (ME) \eqref{eq:numerical_approximation_MES}, the tanh Euler method (TE) \eqref{eq:numerical_approximation_TES} and the sin Euler method (SE) \eqref{eq:numerical_approximation_SES}. For these examples, the coefficients satisfy Assumption 11, where both the drift and diffusion terms are non-globally Lipschitz.

To approximate the law $\mathcal L(X_{t_k})$ at each time step $t_k$ for $k = 0, 1, \dots, n$ by its empirical distribution, we apply the particle method with the number of particles $N = 100$.

As we do not know the exact solution of the considered examples, the strong convergence with respect to the number of time steps is assessed by comparing two solutions computed on a fine and coarse time grid, respectively, using the same samples of Brownian motion. The reference values of the models are computed based on a Monte Carlo method combined with the approximations and $h_{ref} =2^{-17} $. We also apply Monte Carlo method to compute the root-mean-square error (RMSE) in step size $h = \{2^{- 13}, 2^{- 13.5}, 2^{- 14}, 2^{- 14.5}, 2^{- 15}, 2^{- 15.5}, 2^{- 16}\}$ by
$$
\text{RMSE} := \sqrt{\frac{1}{N} \sum_{i = 1}^{N} \left(X_{T}^{i, N, n} - X_{T}^{i, N, n_h} \right)^2}
$$
at the terminal time $T = 1$, where $i$ refers to the $i$-th particle, and $n = T/h_{ref} = 2^{17}$ and $n_h = \lfloor T/h \rfloor$ are the corresponding time steps in the fine and coarse time grid, respectively.

\begin{figure}[htb]
\begin{minipage}{\linewidth}
\centering
\subcaptionbox{}	
{\includegraphics[width=0.42\linewidth, height=0.28\textheight]{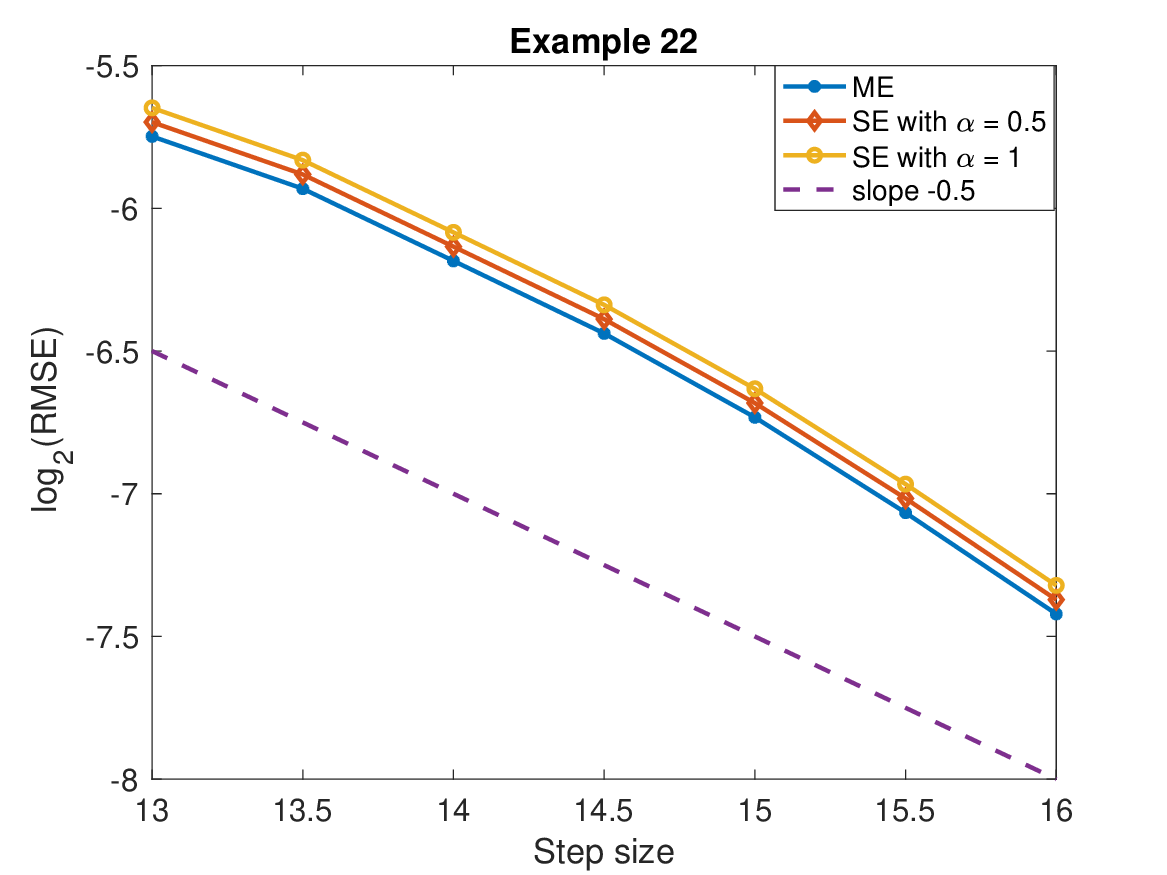}}\quad
\subcaptionbox{}
{\includegraphics[width=0.42\linewidth, height=0.28\textheight]{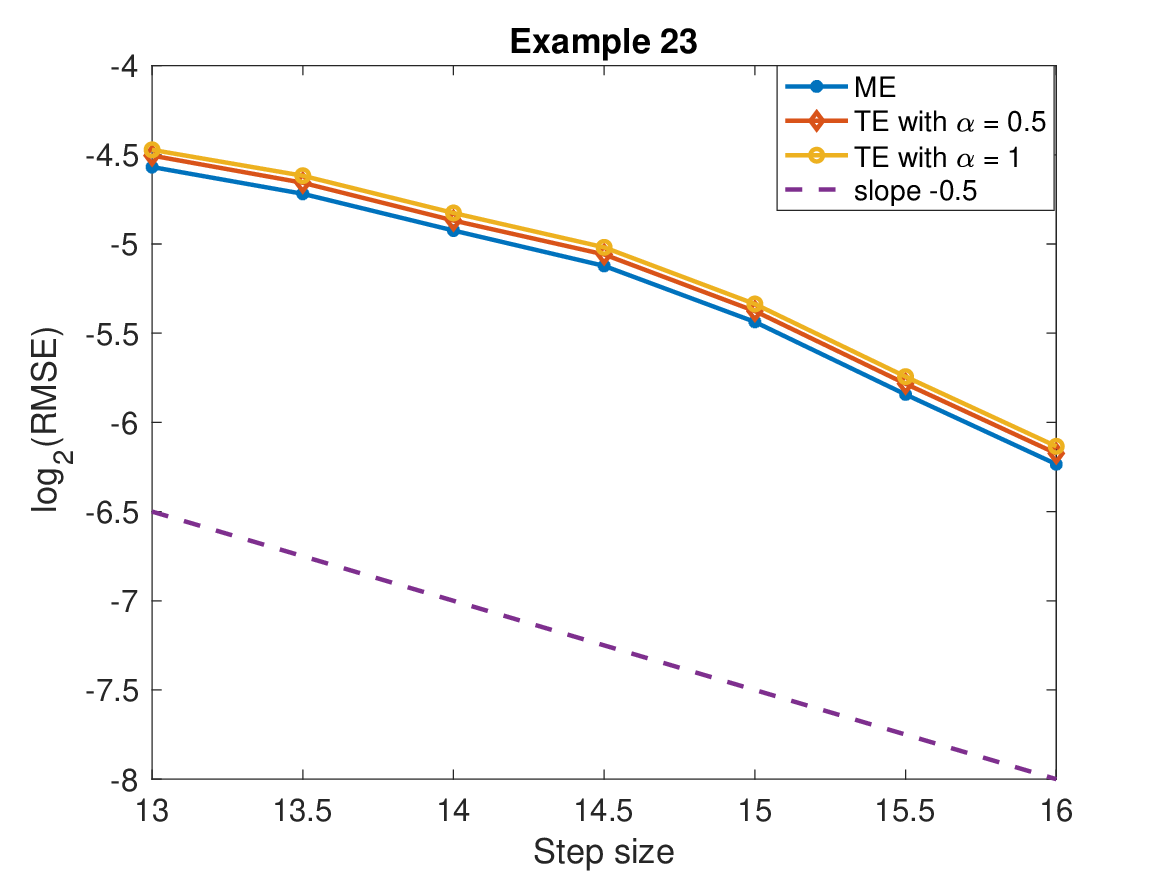}}
\caption{Strong errors for Example \ref{ex:numerical_ex1_nonlinear_diff} and \ref{ex:numerical_ex2_nonlinear_diff}}
\label{fig:strong-error-non-linear-diff}
\end{minipage}
\end{figure}

\begin{exmp}
\label{ex:numerical_ex1_nonlinear_diff}
As the first test model, we consider the following McKean-Vlasov SDE
\begin{equation*}
\label{eq:non_linear_model_1}
\begin{cases}
d X_t = \left( X_t - X_t^3 + c \E[ X_t ] \right) d t + \gamma \left( 1 - X_t^2 \right) d W(t), \quad t \in (0, T], \\
X_0 = x, 
\end{cases}
\end{equation*}
with the parameter values $\gamma = 0.5$, $c = 1$, and $x = 0$. The conditions in Assumption \ref{a:euler_coefficient_b_sigma} are fulfilled with $\rho \geq 1$.
\end{exmp}
For this example, we test the RMSE of three different numerical methods ME \eqref{eq:numerical_approximation_MES}, SE \eqref{eq:numerical_approximation_SES} with $\alpha = 1/2$ and SE \eqref{eq:numerical_approximation_SES} with $\alpha = 1$. As expected, the strong convergence rates of ME \eqref{eq:numerical_approximation_MES}, SE \eqref{eq:numerical_approximation_SES} with $\alpha = 1/2$ and SE \eqref{eq:numerical_approximation_SES} with $\alpha = 1$ are close to $1/2$ from Fig. \ref{fig:strong-error-non-linear-diff} (a).
\begin{exmp}
\label{ex:numerical_ex2_nonlinear_diff}
For the second test model, we consider a McKean-Vlasov SDE in which the drift term preserves higher-order growth condition:
\begin{equation*}
\label{eq:non_linear_model_2}
\begin{cases}
d X_t = \left( 1 - X_t^5 + X_t^3 + c \E[ X_t ] \right) d t + \left(\gamma X_t^2 + 1 \right) d W(t), \quad t \in (0, T], \\
X_0 = x, 
\end{cases}
\end{equation*}
with the parameter values $c = 1$, $\gamma = 0.01$ and $x = 0$. It is clear that Assumption \ref{a:euler_coefficient_b_sigma} is 
satisfied if $\rho \geq 2$.
\end{exmp}
In Fig. \ref{fig:strong-error-non-linear-diff} (b), we reveal the RMSE of three different numerical methods ME \eqref{eq:numerical_approximation_MES}, TE \eqref{eq:numerical_approximation_TES} with $\alpha = 1/2$ and TE \eqref{eq:numerical_approximation_TES} with $\alpha = 1$ against the same step sizes on the $\log_2(\cdot)$ scale. The strong error rate is $1/2$ as expected in Theorem \ref{t:convergence_rate_euler}.

The third example features a non-globally Lipschitz drift and a global Lipschitz diffusion with respect to the state. We test the drift-tamed Euler (DTE) \eqref{eq:numerical_approximation_DTE}, the modified Euler method (ME) \eqref{eq:numerical_approximation_MES}, the tanh Euler method (TE) \eqref{eq:numerical_approximation_TES}, the sin Euler method (SE) \eqref{eq:numerical_approximation_SES} and the  fully-tamed Euler method (FTE) \eqref{eq:numerical_approximation_FTE} facilitating a comparative analysis with an alternative numerical approach, the split-step method (SSM) proposed in \cite{chen2022flexible}, \cite{chen2024euler}, and \cite{chen2023wellposedness}.

\begin{exmp}
\label{ex:linear_diff_ini_norm_distri}
We consider the double-well model considered in \cite{chen2024euler}, given by
\begin{equation*}
 \begin{cases}
d X_t = \left( -\frac54 X_t^3 + 3 X_t^2 \E[ X_t ] -3 X_t \E[ X^2_t ]+\E[ X^3_t ]\right) d t + X_t d W(t), \quad t \in (0, T], \\
X_0 \sim \mathcal{N}(\mu,\,\sigma^{2}), 
\end{cases} 
\end{equation*}
where $\mathcal{N}(\mu, \sigma^2)$ is the normal distribution with mean $\mu$ and variance $\sigma^2$. There are three stable states $\{-2,0,2\}$ for this
model \cite{tugaut2013convergence}. 
\end{exmp}

\begin{figure}
\begin{center}
\includegraphics[height=4.15cm,width=6cm]{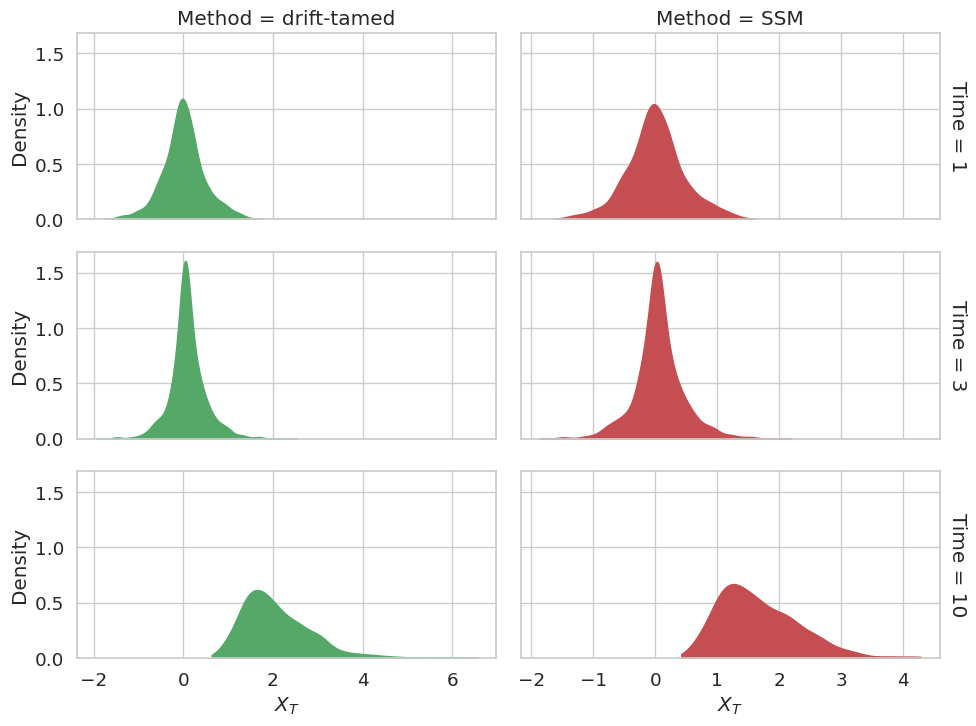}   
\caption{Density with $X_0 \sim \mathcal{N}(0,\,1)$}  
\label{fig:lin-diff-density-drift-ssm}
\end{center} 
\end{figure}

\begin{figure}
\begin{center}
\includegraphics[height=4.15cm, width=6cm]{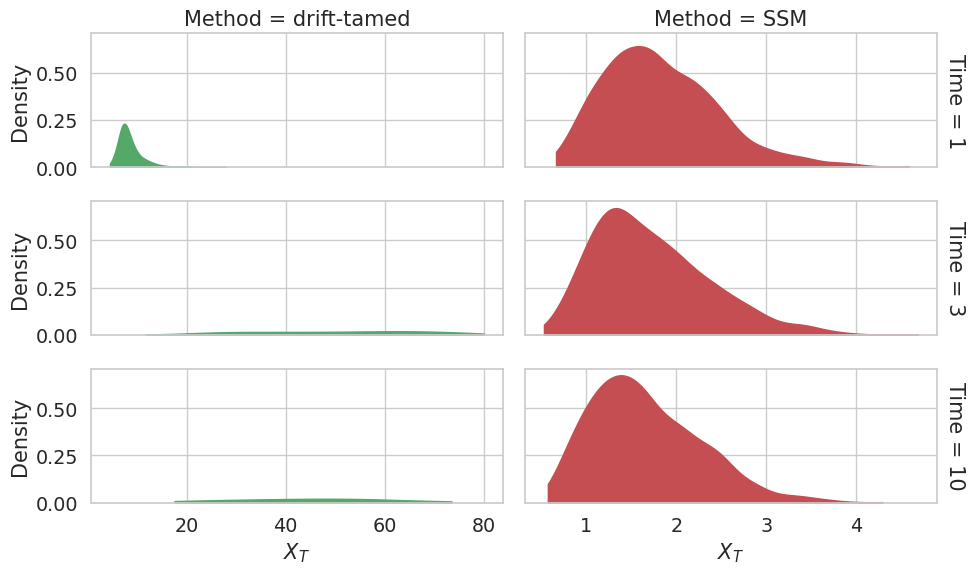} 
\caption{Density with $X_0 \sim \mathcal{N}(3,\,9)$}  
\label{fig:lin-diff-density-drift-ssm-33}
\end{center} 
\end{figure}

\begin{figure}
\begin{center}
\includegraphics[height=4.15cm, width=6cm]{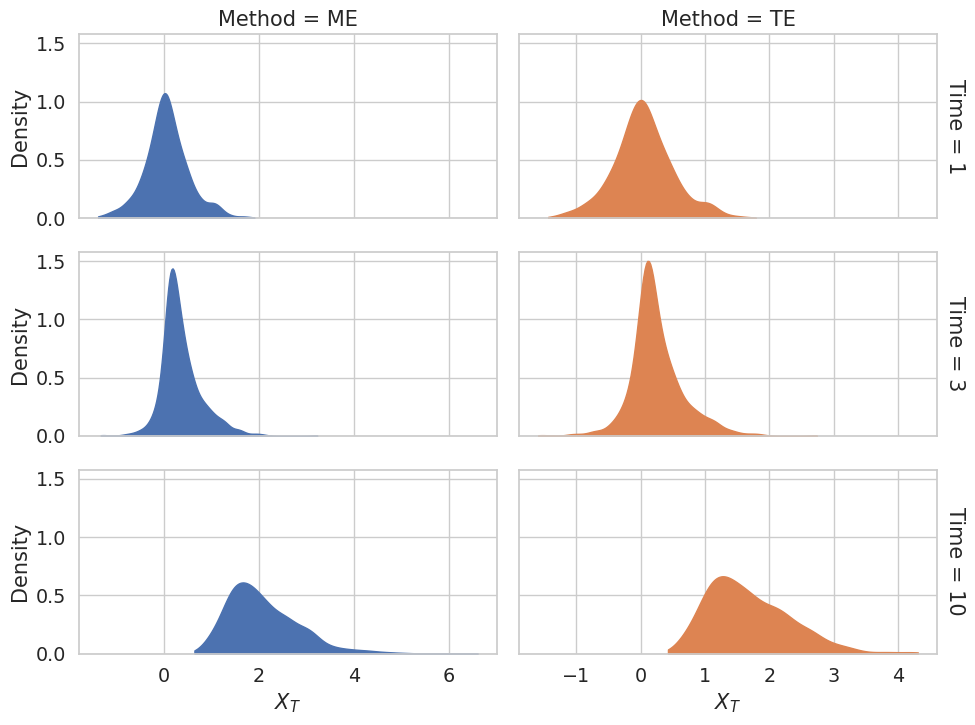}    
\caption{Density with $X_0 \sim \mathcal{N}(0,\,1)$}  
\label{fig:lin-diff-density-ME-TE-01}  
\end{center} 
\end{figure}

\begin{figure}
\begin{center}
\includegraphics[height=4.15cm, width=6cm]{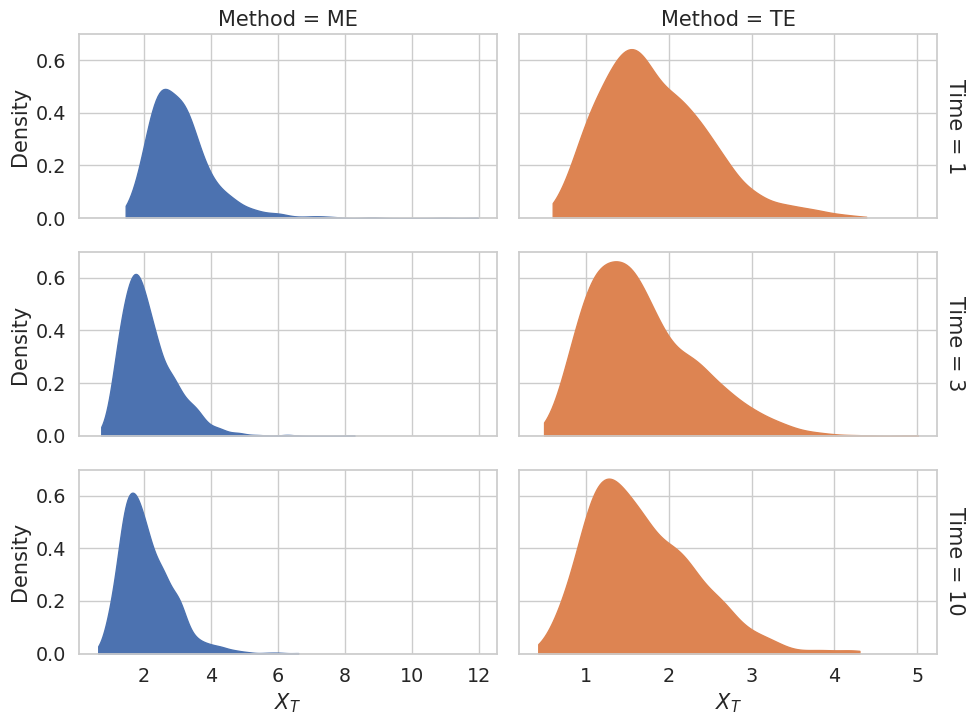}    
\caption{Density with $X_0 \sim \mathcal{N}(3,\,9)$}  
\label{fig:lin-diff-density-ME-TE-33}  
\end{center} 
\end{figure}
\begin{figure}[htb]
\begin{minipage}{\linewidth}
\centering
\subcaptionbox{}	
{\includegraphics[width=0.35\linewidth, height=0.28\textheight]{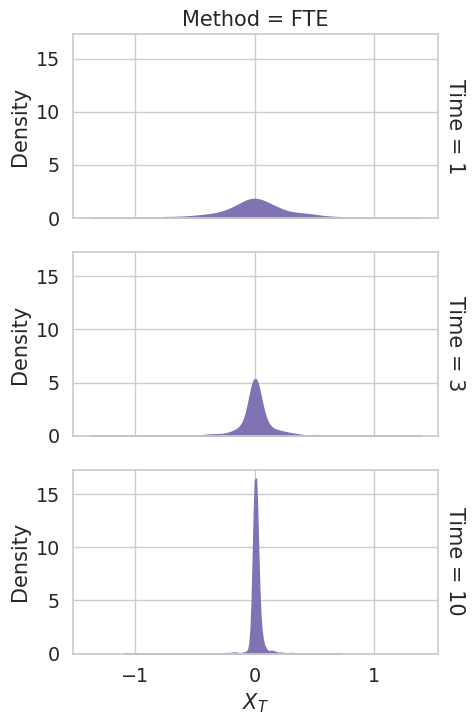}}\quad
\subcaptionbox{}
{\includegraphics[width=0.35\linewidth, height=0.28\textheight]{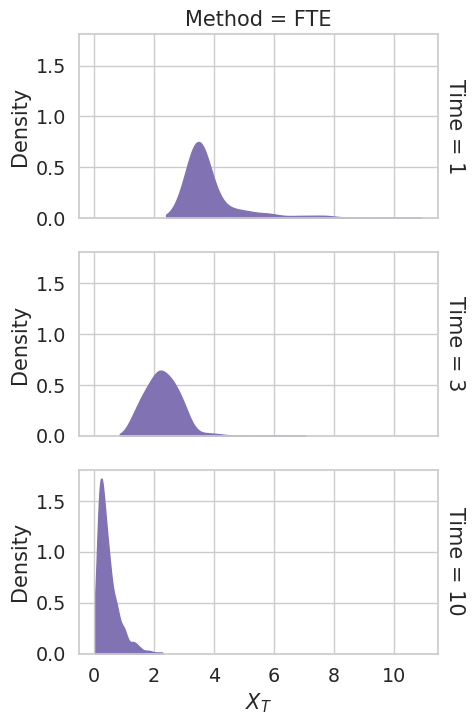}}
\caption{Density with $X_0 \sim \mathcal{N}(0,\,1)$(a) and $X_0 \sim \mathcal{N}(3,\,9)$(b)}
\label{fig:density-fully-tamed}
\end{minipage}
\end{figure}
\begin{figure}[htb]
\begin{minipage}{\linewidth}
\centering
\subcaptionbox{}	
{\includegraphics[width=0.35\linewidth, height=0.30\textheight]{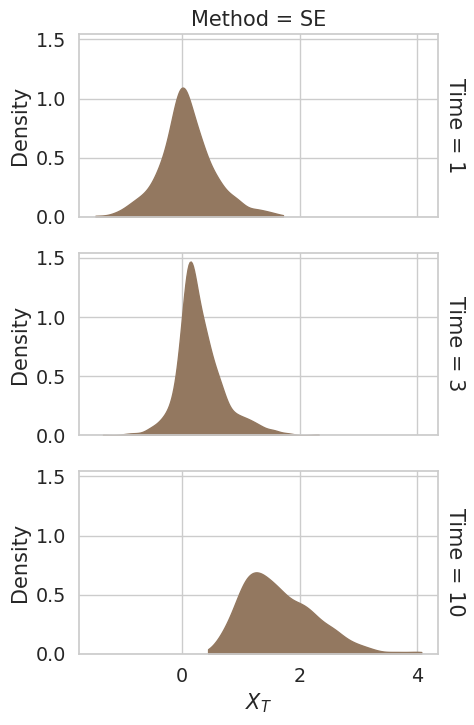}}\quad
\subcaptionbox{}
{\includegraphics[width=0.35\linewidth, height=0.30\textheight]{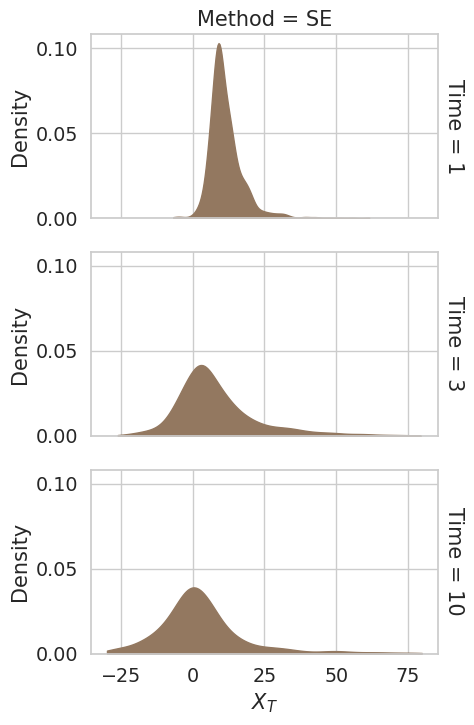}}
\caption{Density with $X_0 \sim \mathcal{N}(0,\,1)$(a) and $X_0 \sim \mathcal{N}(3,\,9)$(b)}
\label{fig:density-SE}
\end{minipage}
\end{figure}

\begin{figure}
\begin{center}
\includegraphics[height=4.0cm,width=6.5cm]{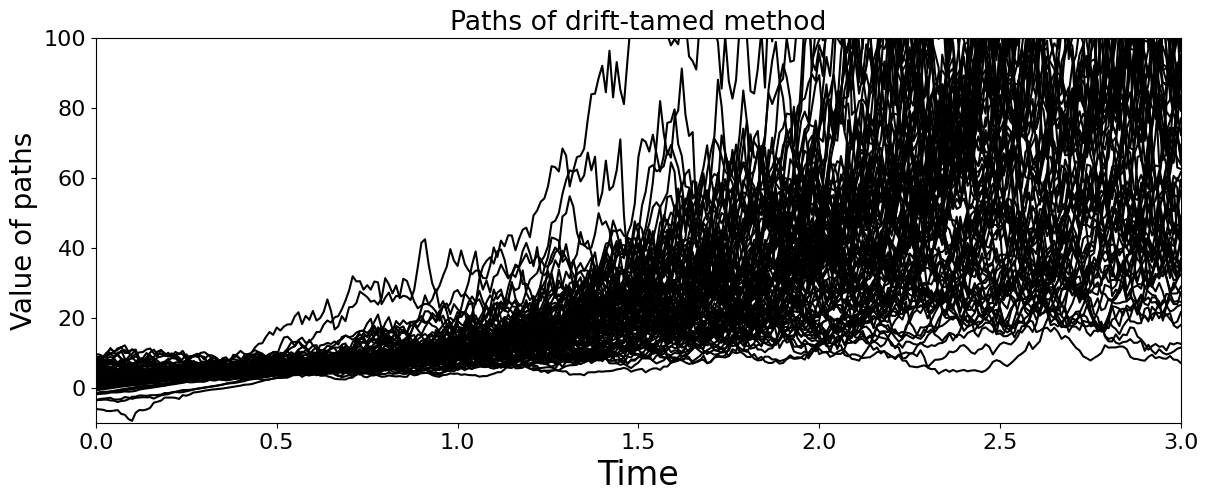}    
\caption{Paths of drift-tamed method for Example \ref{ex:linear_diff_ini_norm_distri}} 
\label{fig:paths-val-drift-tamed}
\end{center}
\end{figure}

\begin{figure}
\begin{center}
\includegraphics[height=4.0cm,width=6cm]{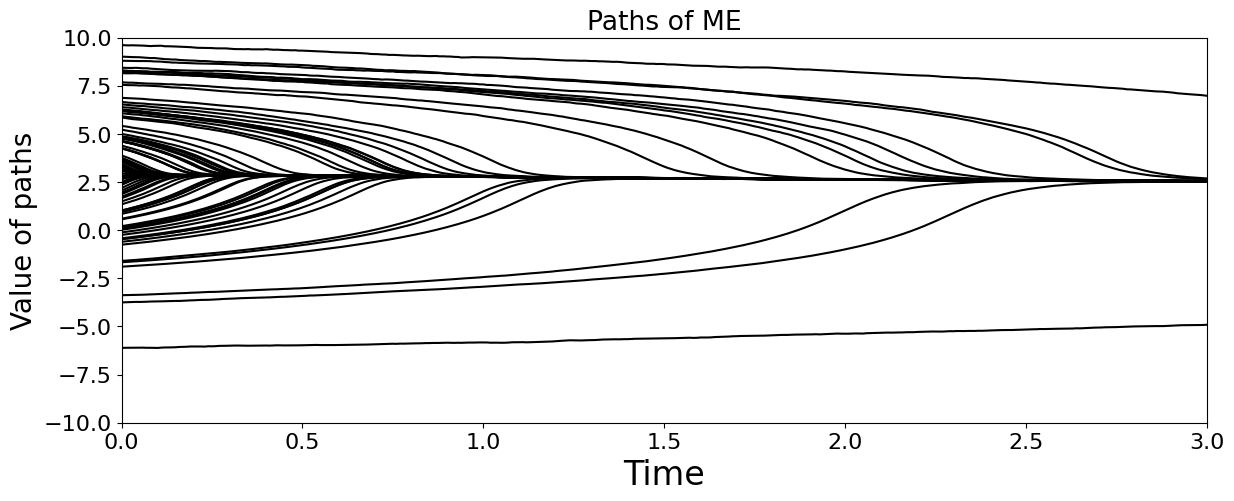}    
\caption{Paths of modified Euler method for Example \ref{ex:linear_diff_ini_norm_distri}}  
\label{fig:paths-val-ME}                   
\end{center}                               
\end{figure}

\begin{figure}
\begin{center}
\includegraphics[height=4.0cm,width=6cm]{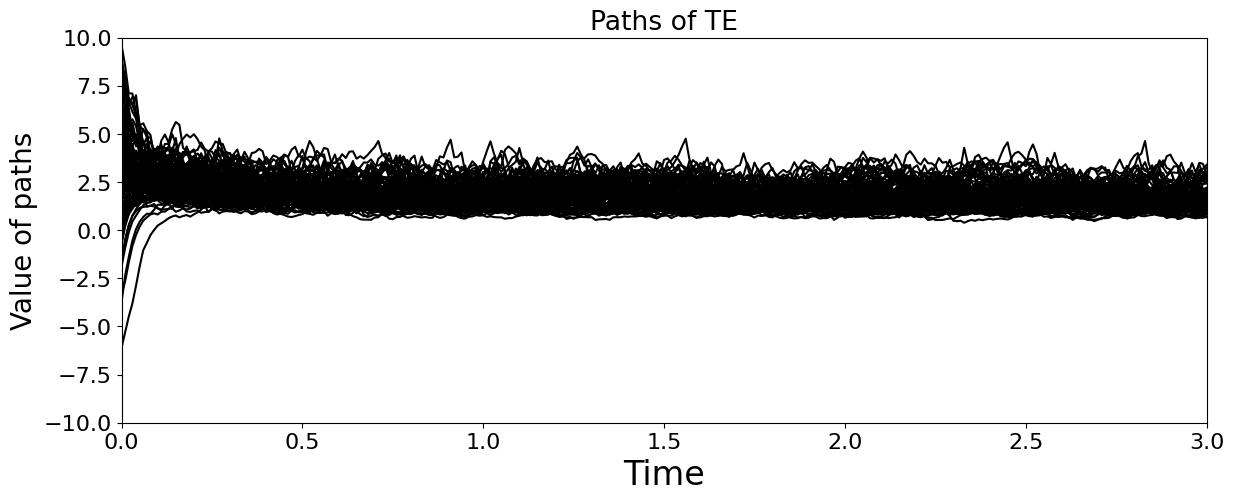}    
\caption{Paths of tanh Euler method for Example \ref{ex:linear_diff_ini_norm_distri}}  
\label{fig:paths-value-TE}      
\end{center}                               
\end{figure}

\begin{figure}
\begin{center}
\includegraphics[height=4.0cm,width=6cm]{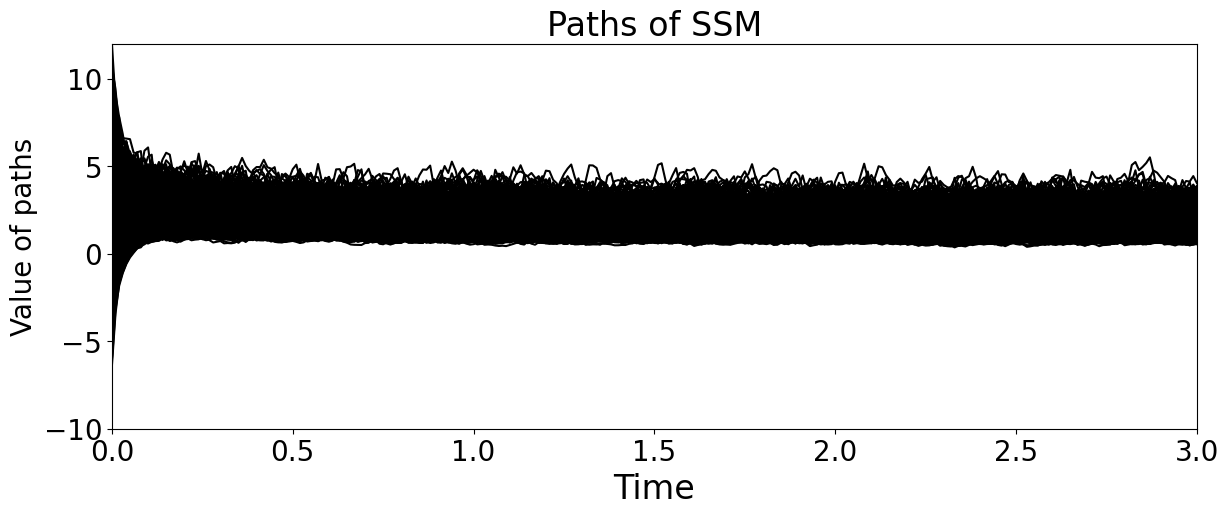}    
\caption{Paths of split-step method for Example \ref{ex:linear_diff_ini_norm_distri}}  
\label{fig:paths-value-SSM}   
\end{center}                               
\end{figure}

\begin{figure}
\begin{center}
\includegraphics[height=4.0cm,width=6cm]{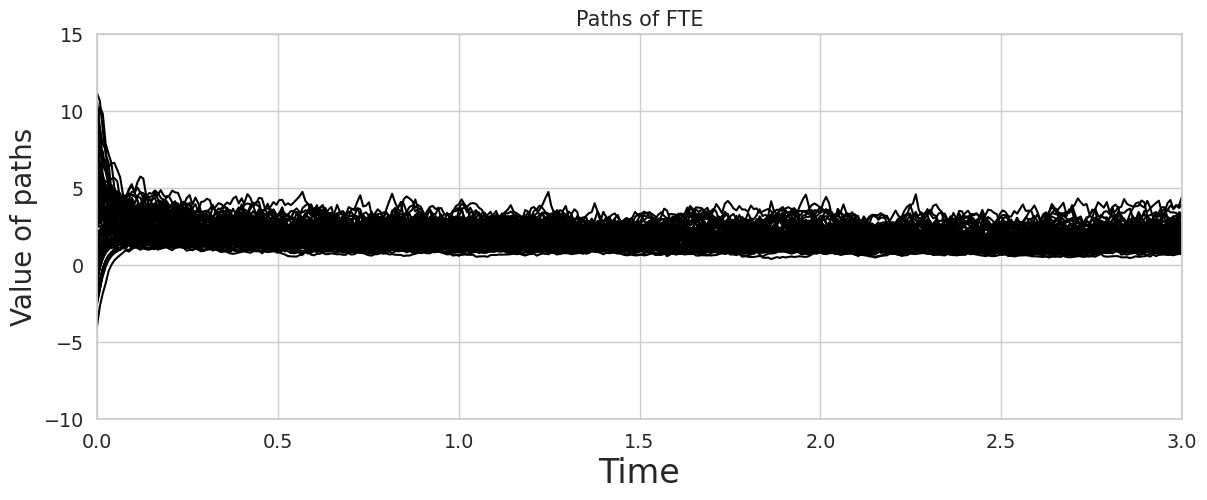}    
\caption{Paths of fully-tamed Euler method for Example \ref{ex:linear_diff_ini_norm_distri}}  
\label{fig:paths-value-FTE}                
\end{center}    
\end{figure}

\begin{figure}
\begin{center}
\includegraphics[height=4.0cm,width=6cm]{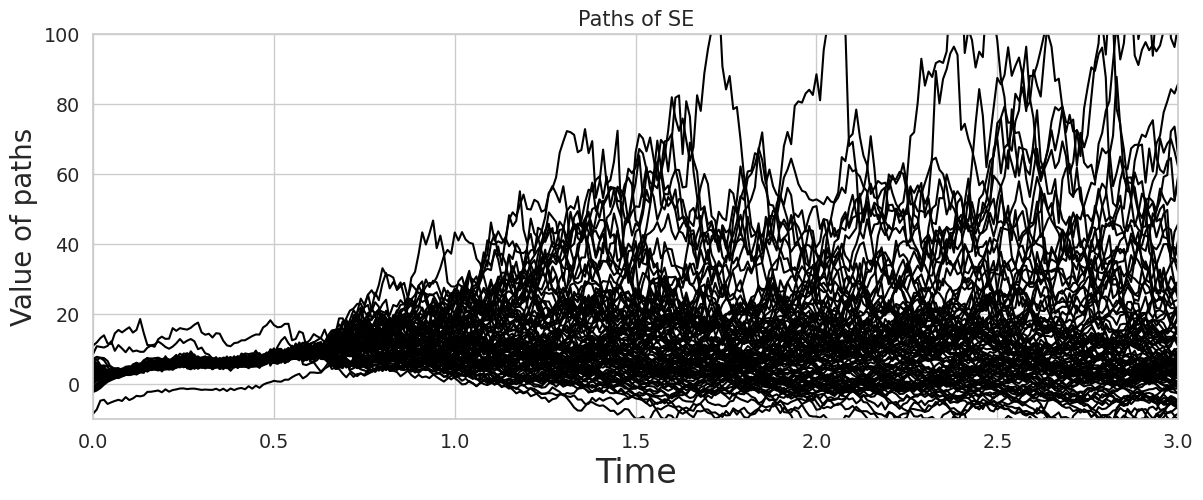}    
\caption{Paths of sin Euler method for Example \ref{ex:linear_diff_ini_norm_distri}}  
\label{fig:paths-value-SE}                 
\end{center}          
\end{figure}

Using the same setup as in \cite{chen2024euler} with $N=1000$, and computing the reference solution with a step size of $h = 10^{-4}$, Figs. \ref{fig:lin-diff-density-drift-ssm}-\ref{fig:density-SE} show density maps for DTE \eqref{eq:numerical_approximation_DTE} with $\lambda=1/2$, ME \eqref{eq:numerical_approximation_MES}, TE \eqref{eq:numerical_approximation_TES} with $\alpha = 1$, SE \eqref{eq:numerical_approximation_SES} with $\alpha = 1$, FTE \eqref{eq:numerical_approximation_FTE} and SSM. These results use a step size of $h=10^{-2}$ and are shown at times $T = 1,3,10$ for two different initial distributions $\mathcal{N}(0,\,1), \mathcal{N}(3,\,9)$. Simulated paths of these methods with the initial distribution $\mathcal{N}(3,\,9)$ are illustrated in Figs. \ref{fig:paths-val-drift-tamed}-\ref{fig:paths-value-SE}.

As noted in \cite{chen2024euler}, for large initial values (i.e., $X_0 \sim \mathcal{N}(3,\,9)$), the DTE \eqref{eq:numerical_approximation_DTE} with 
$\lambda=1/2$ produces unacceptable results (see Fig \ref{fig:lin-diff-density-drift-ssm-33}). This happens because the method \eqref{eq:numerical_approximation_DTE} becomes unstable with the stepsize $h=10^{-2}$, as clearly indicated by Fig. \ref{fig:paths-val-drift-tamed}. 
Similar phenomenon can be also detected for the SE \eqref{eq:numerical_approximation_SES} (see Fig. \ref{fig:density-SE}, \ref{fig:paths-value-SE}).
In contrast, the SSM (\cite{chen2024euler}, \cite{chen2022flexible}), ME \eqref{eq:numerical_approximation_MES}, TE \eqref{eq:numerical_approximation_TES} with $\alpha = 1$ and FTE \eqref{eq:numerical_approximation_FTE} produce acceptable results in the same setting.
Interestingly, by decreasing the stepsize to e.g.,
$h = 0.004$, the DTE \eqref{eq:numerical_approximation_DTE} with 
$\lambda=1/2$ and SE \eqref{eq:numerical_approximation_SES} with $\alpha = 1$ can be then stable and give acceptable approximations.

In terms of the density maps depicted in Figs. \ref{fig:lin-diff-density-drift-ssm}-\ref{fig:density-SE}, one can clearly observe that, the TE \eqref{eq:numerical_approximation_TES} performs most closely to the SSM and gives more reliable approximations than the other explicit methods.

To sum up,
as an implicit method, the SSM method has better stability properties than the other explicit methods and can thus always produce reliable approximations even when treating  MV-SDEs with relatively large initial values and relatively large stepsize $h$. However, some explicit methods, such as the DTE \eqref{eq:numerical_approximation_DTE} and SE \eqref{eq:numerical_approximation_SES} with $\alpha = 1$, might be sensitive to the stepsize selection. Moreover, the above numerical results demonstrate that the TE \eqref{eq:numerical_approximation_TES} performs better than the other explicit methods.

\newpage

\begin{ack}                 
This work was supported by NSF grant of China (No.12471394, 12071488, 12371417) and the Postdoctoral Fellowship Program of CPSF under Grant Number GZC2024205.

We are grateful to four anonymous refrees and the associate editor for providing valuable suggestions.
\end{ack}

\bibliographystyle{plain}        
\bibliography{autosam}           


\appendix
\section{Proof of Lemma \ref{l:moment_bound_euler_discrete_time} in Section \ref{sub:bound-moment-modi-method}}
\label{appen:proof-moment-bound-modi-method}
\textit{{Proof of} Lemma \ref{l:moment_bound_euler_discrete_time}}.
In the following, we use $K$ to denote the generic constant which is independent of $n$ and $h$. Let $\mathcal R > 0$ be sufficiently large and define a sequence of decreasing subevents
$$
\Omega_{\mathcal R, k} = \left\{\omega \in \Omega : \left\vert X_{j}^{i, N, n}(\omega)\right\vert \leq \mathcal R, \, j = 0, 1, \dots, k \right\}
$$
for $k = 0,1, \dots, n$. We denote the complement of $\Omega_{\mathcal R, k}$ by $\Omega_{\mathcal R, k}^{c}$.

Firstly, we show that the boundedness of the moment is valid within a family of appropriate subevents $\left\{\Omega_{\mathcal R, k} \right\}_{k \in \{0,1, \dots, n\}}$. Note that 
\begin{equation}
\label{eq:discom_sub_event_scheme}
\begin{aligned}
& \mathbb{E} \left[\1_{\Omega_{\mathcal{R}, k+1}} \left\vert X_{k+1}^{i, N, n} \right\vert^{2 \bar p} \right] \leq \mathbb{E} \left[\1_{\Omega_{\mathcal R, k}} \left|X_{k+1}^{i, N, n}\right|^{2 \bar p}\right]  \\
= \ & \mathbb{E} \left[\1_{\Omega_{\mathcal R, k}} \left|X_{k+1}^{i, N, n} - X_{k}^{i, N, n} + X_{k}^{i, N, n} \right|^{2 \bar p} \right] \\
\leq \ & \mathbb{E} \left[\1_{\Omega_{\mathcal R, k}} \left|X_{k}^{i, N, n} \right|^{2 \bar p} \right] \\
& \hspace{0.2in} + \mathbb{E} \left[\1_ {\Omega_{\mathcal R ,k}} \left| X_{k}^{i, N, n} \right|^{2 \bar p - 2} \left( 2 \bar p \left \langle X_{k}^{i, N, n}, X_{k+1}^{i, N, n} - X_{k}^{i, N, n} \right\rangle + \bar p \left(2 \bar p - 1 \right) \left|X_{k+1}^{i, N, n} - X_{k}^{i, N, n} \right|^{2} \right) \right] \\
& \hspace{0.2in} + K \sum_{l = 3}^{2 \bar p} \mathbb{E} \left[\1_{\Omega_{\mathcal R, k}} \left|X_{k}^{i, N, n}\right|^{2 \bar p - l} \left|X_{k+1}^{i, N, n} - X_{k}^{i, N, n} \right|^{l} \right] \\
:= \ & \mathbb{E} \left[\1_{\Omega_{\mathcal R, k}} \left|X_{k}^{i, N, n} \right|^{2 \bar p} \right] + I_{1} + I_{2}.
\end{aligned}
\end{equation}
For $I_1$, according to \eqref{eq:modified_euler_scheme} and $\Delta W_{r}^{i} \left(t_{k} \right)$ is independent with $X_{k}^{i, N, n}$ for all $r = 1,2, \dots, m$, one can derive
\begin{equation*}
\begin{aligned}
I_{1} & =  2 \bar p \mathbb{E} \left[\1_{\Omega_{\mathcal R, k}} \left|X_{k}^{i, N, n} \right|^{2 \bar p - 2} \left\langle X_{k}^{i, N, n}, \mathcal T_{1} \left(b \left(t_{k}, X_{k}^{i, N, n}, \mu_{t_k}^{X, N, n} \right), h \right) h -  b \left(t_{k}, X_{k}^{i, N, n}, \mu_{t_k}^{X, N, n} \right) h \right\rangle \right]  \\
& \hspace{0.3in} + 2 \bar p \mathbb{E} \left[\1_{\Omega_{\mathcal R, k}} \left|X_{k}^{i, N, n} \right|^{2 \bar p - 2} \left\langle X_{k}^{i, N, n}, b \left(t_{k}, X_{k}^{i, N, n}, \mu_{t_k}^{X, N, n} \right) h \right\rangle \right]  \\
& \hspace{0.3in} + \bar p \left(2 \bar p - 1 \right) \mathbb{E} \left[\1_{\Omega_{\mathcal R, k}} \left|X_{k}^{i, N, n}\right|^{2 \bar p - 2} \left| \mathcal T_{1} \left( b \left(t_{k}, X_{k}^{i, N, n}, \mu_{t_k}^{X, N, n} \right), h \right) h \right|^2 \right] \\
& \hspace{0.3in} + \bar p \left(2 \bar p - 1 \right) \mathbb{E} \bigg[\1_{\Omega_{\mathcal R, k}} \left|X_{k}^{i, N, n}\right|^{2 \bar p - 2} \bigg|
\sum_{r = 1}^{m} \mathcal T_{2} \left( \sigma_{r} \left( t_{k}, X_{k}^{i, N, n}, \mu_{t_{k}}^{X, N, n} \right), h \right) \Delta W_{r}^{i} \left( t_{k} \right) \bigg|^{2} \bigg].
\end{aligned}
\end{equation*}
Using the Schwartz inequality, assumptions (\textit{A2}), (\textit{H1}), (\textit{H2}) and the growth condition for
the coefficient $b$ in \eqref{eq:growth_condition_b}, we have
\begin{equation*}
\begin{aligned}
I_{1} 
& \leq K h^{1 + r_{1}} \mathbb{E} \left[\1_{\Omega_{\mathcal R, k}} \left|X_{k}^{i, N, n} \right|^{2 \bar p - 1} \left( 1 + \left|X_{k}^{i, N, n} \right|^{2 \rho + 1 } + \mathcal W_{2} \left(\mu_{t_{k}}^{X, N, n}, \delta_{0}\right)\right)^{r_{2}}\right]  \\
& \hspace{0.3in} + K h^{2} \mathbb{E} \left[\1_{\Omega_{\mathcal R, k}} \left|X_{k}^{i, N, n} \right|^{2 \bar p - 2} \left(1+ \left|X_{k}^{i, N, n} \right|^{2 \rho + 1} + \mathcal W_{2} \left(\mu_{t_{k}}^{X, N, n}, \delta_{0} \right) \right)^{2} \right]  \\
& \hspace{0.3in} + K h \mathbb{E} \left[\1_{\Omega_{\mathcal R, k}} \left|X_{k}^{i, N, n} \right|^{2 \bar p - 2} \left(1 + \left|X_{k}^{i, N, n}\right|^{2} + \mathcal W_{2}^{2} \left(\mu_{t_{k}}^{X, N, n}, \delta_{0} \right) \right) \right].
\end{aligned}
\end{equation*}
Note that, by Lemma 2.3 of \cite{DST19},
\begin{equation}
\label{eq:w2_definition}
\mathcal W_{2}^{2} \left(\mu_{t_k}^{X, N, n}, \delta_{0} \right) = \frac{1}{N} \sum_{i=1}^{N} \left|X_{k}^{i, N, n} \right|^{2},
\end{equation}
and some simplification, the estimation for $I_1$ is given as follows
\begin{equation*}
\begin{aligned}
I_{1} & \leq 
K h + K h^{1 + r_{1}} 
\sup_{i \in \{1,2, \dots, N\}} \mathbb{E} \left[\1_{ \Omega_{\mathcal R , k}} \left|X_{k}^{i, N, n} \right|^{2 \bar p - 1 + r_{2} (2 \rho + 1)} \right]  \\
& \hspace{0.3in} + K h^{2} 
\sup_{i \in \{1,2, \dots, N\}} \mathbb{E} \left[\1_{\Omega_{\mathcal R, k}}  \left|X_{k}^{i, N, n}\right|^{2 \bar p +4 \rho}\right] + K h \sup_{i \in\{1,2, \dots, N\}} \mathbb{E}\left[\1_{\Omega_{\mathcal R, k}} \left|X_{k}^{i, N, n} \right|^{2 \bar p} \right].
\end{aligned}
\end{equation*}
Next, we focus on the estimation of $I_2$. By the assumption (\textit{H1}), one can get
\begin{equation*}
\begin{aligned}
I_{2} 
& \leq
K \sum_{l=3}^{2 \bar p} 
\mathbb{E} \bigg[\1_{\Omega_{\mathcal R, k}} |X_{k}^{ i, N, n}|^{2 \bar p - l} \Big( \left| \mathcal T_{1} \left( b \left(t_{k}, X_{k}^{i, N, n}, \mu_{t_{k}}^{X, N, n} \right), h \right) \right|^{l} h^{l}  \\
& \hspace{0.3in} + \Big|\sum_{r=1}^{m} \mathcal T_{2} \left( \sigma_{r} \left(t_{k}, X_{k}^{i, N, n}, \mu_{t_{k}}^{X, N, n} \right), h \right) \Delta W_{r}^{i} \left(t_{k} \right) \Big|^{l}
\Big) \bigg]  \\
& \leq
K \sum_{l=3}^{2 \bar p} \mathbb{E}
\bigg[\1_{\Omega_{\mathcal R, k}} 
|X_{k}^{i, N, n }|^{2 \bar p - l}
\Big(\left| b\left(t_{k}, X_{k}^{i, N, n}, \mu_{t_{k}}^{X, N, n}\right)\right|^{l} h^{l}  \\ 
& \hspace{0.3in} + h^{\frac{l}{2}} \sum_{r = 1}^{m} \left| \sigma_{r} \left(t_{k}, X_{k}^{i, N, n}, \mu_{t_{k}}^{X, N, n} \right) \right|^l  \Big) \bigg].
\end{aligned}
\end{equation*}
Then, the growth condition for the coefficients in \eqref{eq:growth_condition_b} and \eqref{eq:growth_condition_sigma} implies
\begin{equation*}
\begin{aligned}
I_{2} 
& \leq
K \sum_{l=3}^{2 \bar p} \mathbb{E}\left[\1_{\Omega_{\mathcal R, k}}\left|X_{k}^{i, N, n}\right|^{2 \bar p -l} \left(1+\left|X_{k}^{i, N, n}\right|^{2 \rho+1} + \mathcal W_{2}\left(\mu_{t_{k}}^{X, N, n}, \delta_{0}\right)\right)^{l} h^{l}\right] \\
& \hspace{0.3in}  + K \sum_{l=3}^{2 \bar p} \mathbb{E}\left[\1_{\Omega_{\mathcal R, k}}\left|X_{k}^{i, N, n}\right|^{2 \bar p -l} \left(1+\left|X_{k}^{i, N, n} \right|^{\rho+1} + \mathcal W_{2}\left(\mu_{t_k}^{X, N, n}, \delta_{0} \right)\right)^{l} h^{\frac{l}{2}}\right] \\
& \leq
K \sum_{l=3}^{2 \bar p} \mathbb{E}
\left[\1_{ \Omega_{\mathcal R, k}} \left|X_{k}^{ i, N, n} \right|^{2 \bar p - l} \left(h^{l}+h^{\frac{l}{2}} + h^{l}\left|X_{k}^{i, N, n} \right|^{(2 \rho+1) l} + h^{l} \mathcal W_{2}^{l} \left( \mu_{t_k}^{X, N, n}, \delta_{0} \right) 
\right.\right. \\
& \hspace{0.3in} \left.\left.+
 h^{\frac{l}{2}}\left|X_{k}^{i, N, n}\right|^{(\rho+1) l} + h^{\frac{l}{2}} \mathcal W_{2}^{l}\left(\mu_{t_{k}}^{X, N, n}, \delta_{0}\right)\right)\right].
 \end{aligned}
\end{equation*}
It follows from \eqref{eq:w2_definition} that
\begin{equation*}
I_{2} \leq
K h^{\frac{3}{2}} +
K \sum_{l=3}^{2 \bar p} 
\sup_{i \in\{1,2, \dots, N\}} \mathbb{E}\left[\1_{\Omega_{\mathcal R, k}} h^{l}\left|X_{k}^{i, N, n}\right|^{2 \rho l+2 \bar p}\right] + K \sum_{l=3}^{2 \bar p} 
\sup_{i \in\{1,2, \dots, N\}} \mathbb{E}\left[\1_{\Omega_{
\mathcal R, k}} h^{\frac{l}{2}}\left|X_{k}^{i, N, n}\right|^{\rho l+2 \bar p }\right].
\end{equation*}
Combining the above results, we have
\begin{equation*}
\begin{aligned}
& \sup_{i \in \{1,2, \dots, N\}} \mathbb{E} \left[\1_{\Omega_{\mathcal R, k}} \left|X_{k+1}^{i, N, n} \right|^{2 \bar p} \right]  \\
\leq \ &
K h + K h^{2} \sup_{i \in \{1,2, \dots, N\}} \mathbb{E} \left[\1_{\Omega_{\mathcal R, k}} \left|X_{k}^{i, N, n} \right|^{2 \bar p + 4 \rho} \right] + K h^{1 + r_{1}} \sup_{i \in \{1,2, \dots, N\}} \mathbb{E} \left[\1_{\Omega_{
\mathcal R, k}} \left|X_{k}^{i, N, n} \right|^{2 \bar p - 1 + r_{2} (2 \rho+1)} \right] \\
& \hspace{0.2in} + (1 + K h) \sup_{i \in \{1,2, \dots, N\}} \mathbb{E} \left[\1_{\Omega_{\mathcal R, k}} \left|X_{k}^{i, N, n} \right|^{2 \bar p} \right] + K \sum_{l=3}^{2 \bar p} 
\sup_{i \in \{1,2, \dots, N\}} \mathbb{E} \left[\1_{\Omega_{\mathcal R, k}} h^{l} \left|X_{k}^{i, N, n} \right|^{2 \rho l + 2 \bar p} \right] \\
& \hspace{0.2in} + 
K \sum_{l=3}^{2 \bar p} 
\sup_{i \in \{1,2, \dots, N\}} \mathbb{E} \left[\1_{\Omega_{\mathcal R, k}} h^{\frac{l}{2}} \left|X_{k}^{i, N, n} \right|^{\rho l + 2 \bar p} \right].
\end{aligned}
\end{equation*}
Choosing $\mathcal{R} = \mathcal{R}(h)=h^{-1/\mathcal G(\rho, r_1, r_2)}$ with 
$\mathcal G := \mathcal G(\rho, r_1, r_2)$ given as \eqref{eq:constant_G}, for all $l = 3, 4, \dots, 2 \bar{p}$, we have the following inequalies
\begin{equation*}
\begin{aligned}
\1_{\Omega_{\mathcal R, k}}
\left|X_{k}^{i, N, n}\right|^{(2 \rho+1) r_{2}+2 \bar p-1} h^{r_{1}} & = \1_{\Omega_{\mathcal R, k}}\left|X_{k}^{i, N, n}\right|^{2 \bar p}\left(\1_{\Omega_{\mathcal R, k}}\left|X_{k}^{i, N, n}\right|^{\frac{(2 \rho+1) r_{2}-1}{r_{1}}} h\right)^{r_{1}} \leq  K \1_{\Omega_{\mathcal R, k}}\left|X_{k}^{i, N, n}\right|^{2 \bar p}, \\
\1_{\Omega_{\mathcal R, k}}\left|X_{k}^{i, N, n}\right|^{4 \rho+2 \bar p} h & = \1_{\Omega_{\mathcal R, k}} \left|X_{k}^{i, N, n}\right|^{2 \bar p}\left(\1_{\Omega_{\mathcal R, k}}\left|X_{k}^{i, N, n}\right|^{4 \rho} h \right) \leq K \1_{\Omega_{\mathcal R, k}}\left|X_{k}^{i, N, n}\right|^{2 \bar p}, \\
\1_{\Omega_{\mathcal R, k}} \left|X_{k}^{i, N, n}\right|^{2 \rho l+2 \bar p} h^{l-1} & = \1_{\Omega_{\mathcal R, k}}\left|X_{k}^{i, N, n}\right|^{2 \bar p} \left(\1_{\Omega_{\mathcal R, k}}\left|X_{k}^{i, N, n}\right|^{\frac{2 \rho l}{l-1}} h \right)^{l-1} \leq K \1_{\Omega_{\mathcal R, k}}\left|X_{k}^{i, N, n}\right|^{2 \bar p},
\end{aligned}
\end{equation*}
and
$$
\1_{\Omega_{\mathcal R, k}}
\left|X_{k}^{i, N, n}\right|^{\rho l+2 \bar p} h^{\frac{l}{2}-1} = \1_{\Omega_{\mathcal R, k}}\left|X_{k}^{i, N, n}\right|^{2 \bar p}\left(\1_{\Omega_{\mathcal R, k}}\left|X_{k}^{i, N, n}\right|^{\frac{2 \rho l}{l-2}} h\right)^{\frac{l-2}{2}}
\leq K \1_{\Omega_{\mathcal R, k}}\left|X_{k}^{i, N, n}\right|^{2 \bar p},
$$
where $ K $ is independent of $ h $. Thus, we obtain
\begin{equation}
\label{eq:sub_k_numerical_mom_bound}
\sup_{i \in\{1,2, \dots, N\}} \mathbb{E} \left[\1_{\Omega_{
\mathcal R, k}} \left|X_{k+1}^{i, N, n} \right|^{2 \bar p} \right]
\leq K h + (1+K h) \sup_{i \in\{1,2, \dots, N\}} \mathbb{E}\left[\1_{\Omega_{\mathcal R, k}}\left|X_{k}^{i, N, n}\right|^{2 \bar p} \right],
\end{equation}
which implies that
\begin{equation*}
\sup_{i \in\{1,2, \dots, N\} } \mathbb{E} \left[\1_{\Omega_{\mathcal R, k+1}}\left|X_{k+1}^{i, N, n}\right|^{2 \bar p} \right] \leq K h +(1+K h)
\sup_{i \in\{1,2, \dots, N\}}  \mathbb{E}\left[\1_{\Omega_{\mathcal R, k}} \left|X_{k}^{i, N, n}\right|^{2 \bar p} \right].
\end{equation*}
Therefore, by induction or Gr\"onwall inequality in discrete time case, we have
\begin{equation}
\label{eq:subev_numerical_mom_bound}
\begin{aligned}
\sup_{i \in\{1,2, \dots, N\}} \mathbb{E}\left[\1_{\Omega_{\mathcal R, k}} \left|X_{k}^{i, N, n}\right|^{2 \bar p} \right] &
\leq (1+Kh)^{k}\left(1+\mathbb{E}
\left[\left|X_{0}\right|^{2 \bar p} \right]\right)
\leq e^{K h k} \left(1+\mathbb{E}\left[\left|X_{0}\right|^{2 \bar p}\right]\right) \\
& \leq K \left(1+\mathbb{E}\left[\left|X_{0}\right|^{2 \bar p} \right]\right),
\end{aligned}
\end{equation}
where $ K $ is a generic constant that is independent of $h$ and $n$. It remains to estimate $\mathbb{E}[\1_{\Omega_{\mathcal R, k}^{c}} |X_{k}^{i, N, n} |^{2 p}]$. It follows from \eqref{eq:modified_euler_scheme} and (\textit{H1}) that
\begin{equation}
\label{eq:numerical_gene_estimate}
\begin{aligned}
& \left|X_{k+1}^{i, N, n}\right| \\
\leq \ &
\left|X_{k}^{i, N, n} \right| + \left|\mathcal T_{1}\left(b\left(t_{k}, X_{k}^{i, N, n}, \mu_{t_{k}}^{i, N, n} \right), h\right) h \right|
+ \sum_{r=1}^{m} \left|\mathcal T_{2} \left(\sigma_{r}
\left(t_{k}, X_{k}^{i, N, n}, \mu_{t_{k}}^{i, N, n}\right), h \right) \Delta W_{r}^{i} \left(t_{k} \right) \right| \\
\leq \ &
\left|X_{k}^{i, N, n} \right|
+ L h^{-1} + \sum_{r=1}^{m} 
L h^{-\frac{3}{2}} \left|W_{r}^{i}\left(t_{k+1} \right) - W_{r}^{i}\left(t_{k} \right) \right| \\
\leq \ &
\left|X_{0}^{i} \right| + L(k+1) h^{-1} + \sum_{j=0}^{k} \sum_{r=1}^{m} L h^{-\frac{3}{2}} \left|W_{r}^{i}\left(t_{j+1} \right) - W_{r}^{i} \left(t_{j} \right) \right|
\end{aligned}
\end{equation}
by induction. Note that
\begin{equation}
\label{eq:indicator_function_cal}
\1_{\Omega_{\mathcal R, k } ^{c}} = 1 - 
\1_{\Omega_{\mathcal R, k}} = 1 -
\1_{\Omega_{\mathcal R, k-1}}
\1_{\left|X_{k}^{i, N, n}\right| \leq \mathcal R} = \sum_{j=0}^{k} \1_{\Omega_{ \mathcal R, j-1}} \1_{\left|X_{j}^{i, N, n} \right| > \mathcal R},
\end{equation}
where we set $\1_{\Omega_{\mathcal R, -1}}= 1$. Then, applying \eqref{eq:indicator_function_cal}, H\"older's inequality with $\frac{1}{p_1} + \frac{1}{q_1} = 1$ for $q_1 = \frac{2 \bar p}{(4 p + 1) \mathcal G} > 1$ due to $p \leq \frac{2 \bar p - \mathcal G}{2 + 4 \mathcal G}$, and the Markov inequality, we derive that
\begin{equation*}
\begin{aligned}
\mathbb{E}\left[\1_{\Omega_{ \mathcal R, k}^{c}} \left|X_{k}^{i, N, n}\right|^{2 p} \right] & = \sum_{j=0}^{k} \mathbb{E} \left[ 
\left|X_{k}^{i, N, n}\right|^{2 p} 
\1_{\Omega_{\mathcal R, j-1}} \1_{\left|X_{j}^{i, N, n}\right| > \mathcal R} \right] \\
& \leq \sum_{j=0}^{k}\left(\mathbb{E}\left[
\left|X_{k}^{i, N, n}\right|^{2 p p_{1}} \right]\right)^{\frac{1}{p_{1}}} \left(\mathbb{E}\left[\1_{\Omega_{\mathcal R, j-1}} \1_{\left|X_{j}^{i, N, n} \right| > \mathcal R} \right] \right)^{\frac{1}{q_{1}}} \\
& =
\left(\mathbb{E}\left[\left|X_{k}^{i, N, n}\right|^{2 p p_{1}}\right]\right)^{\frac{1}{p_{1}}} \sum_{j=0}^{k}\left(\mathbb{P}\left(
\1_{\Omega_{\mathcal R, j-1}}\left|X_{j}^{i, N, n}\right| > \mathcal R \right)\right)^{\frac{1}{q_{1}}} \\
& \leq \left(\mathbb{E}\left[\left|X_{k}^{i, N, n} \right|^{2 p p_{1}}\right]\right)^{\frac{1}{p_{1}}} \sum_{j=0}^{k} \frac{\left(\mathbb{E}\left[\1_{\Omega_{
\mathcal R, j-1}}\left|X_{j}^{i, N, n} \right|^{2 \bar p}\right]\right)^{\frac{1}{q_{1}}}}{\mathcal R^{\frac{2 \bar p}{q_{1}}}}.
\end{aligned}
\end{equation*}
Note that, $p \leq \frac{2 \bar p - \mathcal G}{2 + 4 \mathcal G}$ implies that $\bar p \geq p p_{1}$. Then, H\"older's inequality and \eqref{eq:numerical_gene_estimate} give
\begin{equation*}
\begin{aligned}
& \left(\mathbb{E}\left[\left|X_{k}^{i, N, n}\right|^{2 p p_{1}} \right]\right)^{\frac{1}{p_{1}}} \leq  \left(\mathbb{E}\left[\left|X_{k}^{i, N, n}\right|^{2 \bar p} \right]\right)^{\frac{p}{\bar p}} \\
\leq \ & K \left(\mathbb{E} \left[\left|X_{0}^{i}\right|^{2 \bar p} 
\right] + \left(k h^{-1}\right)^{2 \bar p} + h^{-3\bar p} \mathbb{E} \left[ \left( \sum_{j=0}^{k-1} \sum_{r=1}^{m} \left| W_{r}^{i}(t_{j+1})-W_{r}^{i}(t_{j}) \right| \right)^{2\bar p} \right] \right)^{\frac{p}{\bar p}} \\
\leq \ &
K h^{-4 p} + K \left(1+\mathbb{E}\left[\left|X_{0}\right|^{2 \bar p}\right]\right)^{\frac{p}{\bar p}}.
\end{aligned}
\end{equation*}
According to the inequalities  \eqref{eq:sub_k_numerical_mom_bound} and \eqref{eq:subev_numerical_mom_bound}, we have
\begin{equation*}
\sup_{i\in\{1,2, \dots, N\}} \mathbb{E}\left[\1_{\Omega_{\mathcal R, k}} \left|X_{k+1}^{i, N, n}\right|^{2 \bar p} \right] \leq
K h+(1+ K h)K \left( 1+\mathbb{E} \left[\left|X_{0}\right|^{2 \bar p} \right]\right) \leq K \left(1 + \mathbb{E}\left[\left|X_{0}\right|^{2 \bar p } \right]\right).
\end{equation*}
Recall that $\mathcal R = h^{-\frac{1}{\mathcal G}}$, then $\mathcal R^{\frac{2 \bar p}{q_{1}}} = h^{-\frac{2 \bar p}{q_{1} \mathcal G}} $. Moreover, the inequality $\frac{2 \bar p}{q_{1} \mathcal G} \geq 4p + 1$ holds as $q_1 = \frac{2 \bar p}{(4 p + 1) \mathcal G}$. Thus,
\begin{equation}
\label{eq:complement_estimation}
\begin{aligned}
& \mathbb{E} \left[\1_{\Omega_{
\mathcal R, k}^{c}}\left|X_{k}^{i, N, n} \right|^{2 p}\right] \\
\leq \ & \left(\mathbb{E}\left[\left|X_{k}^{i, N, n} \right|^{2 p p_{1}}\right]\right)^{\frac{1}{p_{1}}} \sum_{j=0}^{k} \frac{\left(\mathbb{E}\left[\1_{\Omega_{
\mathcal R, j-1}}\left|X_{j}^{i, N, n} \right|^{2 \bar p}\right]\right)^{\frac{1}{q_{1}}}}{\mathcal R^{\frac{2 \bar p}{q_{1}}}} \\
\leq \ &
\left( K h^{-4 p} + K \left(1+\mathbb{E}\left[\left|X_{0}\right|^{2\bar p}\right]\right)^{\frac{p}{\bar p}}\right) 
(1+k) h^{\frac{2 \bar p}{q_{1} \mathcal G}}  K \left(1+\mathbb{E}\left[\left|X_{0}\right|^{2 \bar p}\right]\right)^{\frac{1}{q_{1}}} \\
\leq \ & K \left( 1+\mathbb{E}\left[\left|X_{0}\right|^{2 \bar p}\right]\right)^{\frac{1}{q_{1}}+\frac{p}{\bar p}}.
\end{aligned}
\end{equation}
Therefore, by using of H\"older's inequality, \eqref{eq:subev_numerical_mom_bound} and \eqref{eq:complement_estimation}, we obtain the desired result as follows
\begin{equation*}
\begin{aligned}
& \sup_{i \in\{1,2, \dots, N\}} \mathbb{E}\left[\left|X_{k}^{i, N, n } \right|^{2 p}\right] \\
= \ &
\sup_{i \in\{1,2, \dots, N\}} \mathbb{E}\left[\1_{\Omega_{\mathcal R, k} } \left|X_{k}^{i, N, n}\right|^{2 p} \right] + \sup_{i \in\{1,2, \dots, N\}} \mathbb{E} \left[\1_{\Omega_{\mathcal R, k}^{c}} \left|X_{k}^{i, N, n}\right|^{2 p}\right] \\
\leq \ &
\left(\sup_{i \in\{1,2, \dots, N\}} \mathbb{E}\left[\1_{\Omega_{
\mathcal R, k}}\left|X_{k}^{i, N, n} \right|^{2 \bar p} \right]\right)^{\frac{p}{\bar p}} + \sup_{i \in\{1,2, \dots, N\}} \mathbb{E}\left[\1_{\Omega_{\mathcal R, k} 
 ^{c}} \left|X_{k}^{i, N, n}\right|^{2 p}\right] \\
\leq \ & K \left(1+\mathbb{E}\left[\left|X_{0}\right|^{2\bar p}\right]\right)^{\frac{p}{\bar p}} + K\left(1+\mathbb{E}\left[\left|X_{0}\right|^{2 \bar p}\right]\right)^{\frac{1}{q_{1}}+\frac{p}{\bar p}} \\
\leq \ &
K \left( 1 +  \mathbb{E}\left[\left|X_{0}\right|^{2 \beta p}\right] \right),
\end{aligned}
\end{equation*}
where
$\beta = 1 + \frac{\bar p}{p q_{1}} > 1$.
\qed

\end{document}